\documentclass[11pt]{amsart}
\usepackage{amsfonts}
\usepackage{amssymb}
\usepackage{amsmath}
\usepackage{bbm}
\usepackage[dvipsnames]{xcolor}
\usepackage{tikz}
\usetikzlibrary{trees}
\usetikzlibrary{positioning, calc}
\tikzstyle{vertex}=[inner sep=0pt]
\usetikzlibrary{arrows}
\usetikzlibrary{cd}
\tikzset{>=latex}
\tikzset{LA/.style = {line width=#1, -{Straight Barb[length=3pt]}},LA/.default=1.2pt}

\usetikzlibrary{decorations.pathmorphing}
\usetikzlibrary{decorations.pathreplacing}
\usetikzlibrary{decorations.markings}
\usetikzlibrary{decorations.shapes}
\usepackage{verbatim}
\usepackage[font=small,skip=0pt]{caption}

\usepackage{mathtools}

\usepackage{hyperref}
\usepackage{amsrefs}

\newcommand{\stringdiagram}[1]{\[\begin{tikzpicture}[scale=.33, thick]
#1\end{tikzpicture}\]}

\newcommand{\identity}[2]{\draw (#1,#2+1) -- (#1,#2-1);}
\newcommand{\halfidentity}[2]{\draw (#1,#2) -- (#1,#2-1);}
\newcommand{\unit}[2]{\draw (#1,#2) circle [radius=.2]; \draw (#1,#2-.2) -- (#1,#2-1);}

\newcommand{\morphism}[3]{\draw (#1-.4, #2-.6) rectangle (#1+.4,#2+.6); \draw (#1,#2-.6) -- (#1,#2-1); \draw (#1,#2+.6) -- (#1, #2+1); \node at (#1,#2) {$#3$};}
\newcommand{\multiplication}[2] {\draw (#1,#2-.3) -- (#1-1,#2+1); \draw (#1,#2-.3) -- (#1+1,#2+1); \draw (#1,#2-.3) -- (#1,#2-1);}
\newcommand{\slimmultiplication}[2] {\draw (#1,#2-.3) -- (#1-.5,#2+1); \draw (#1,#2-.3) -- (#1+.5,#2+1); \draw (#1,#2-.3) -- (#1,#2-1);}
\newcommand{\widemultiplication}[2] {\draw (#1,#2-.3) -- (#1-1.5,#2+1); \draw (#1,#2-.3) -- (#1+1.5,#2+1); \draw (#1,#2-.3) -- (#1,#2-1);}
\newcommand{\tripleprod}[2] {\draw (#1,#2) -- (#1-1,#2+1); \draw (#1,#2) -- (#1+1,#2+1); \draw (#1,#2+1) -- (#1,#2-1);}

\newcommand{\equals}[2]{\node at (#1,#2) {$=$};}
\newcommand{\nattrans}[3]{\node at (#1,#2) {$\xRightarrow{#3}$};}
\newcommand{\leftnattrans}[3]{\node at (#1,#2) {$\xLeftarrow{#3}$};}

\newcommand{\highlight}[4]{\draw[gray, ultra thin] (#1,#2) rectangle (#3,#4);}
\newcommand{\enclose}[4]{\draw[brown, thin] (#1,#2) rectangle (#3,#4);}
\newcommand{\boxA}[3]{\draw[draw=black] (#1,#2) rectangle node {${#3}$} ++(3,2);}
\newcommand{\longboxA}[3]{\draw[draw=black] (#1,#2) rectangle node {${#3}$} ++(7,2);}
\newcommand{\dbox}[2]{\draw[dashed, MidnightBlue] (#1,#2) rectangle node {} ++(4,3);}
\newcommand{\longdbox}[2]{\draw[dashed, MidnightBlue] (#1,#2) rectangle node {} ++(8,3);}
\newcommand{\strDots}[2]{\node at (#1,#2) {$\cdots$};}
\newcommand{\strOverBrace}[3]{\draw [decoration={brace}, decorate] (#1,#2) -- (#1+1.5, #2) node [ label, yshift=2ex, pos=0.5] {#3};}

\newcommand{\cat}{\mathcal{C}}
\DeclareMathOperator{\Ob}{Ob}
\DeclareMathOperator{\Mor}{Mor}
\DeclareMathOperator{\Hom}{Hom}
\newcommand{\catname}[1]{\mathbf{#1}}

\newcommand{\Span}{\catname{Span}}

\newcommand{\vect}{\catname{Vect}}
\newcommand{\Cat}{\catname{Cat}}
\newcommand{\set}{\catname{Set}}

\newcommand{\id}{\mathrm{id}}

\newcommand{\tolabel}[1]{\xrightarrow{#1}}
\newcommand{\fromlabel}[1]{\xleftarrow{#1}}

\newcommand{\op}{\mathrm{op}}

\newcommand{\isoto}{\stackrel{\sim}{\to}}
\newcommand{\tame}{\ell}

\newcommand{\suchthat}{\mid}

\newcommand{\N}{\mathbb{N}}
\newcommand{\kk}{\Bbbk}

\renewcommand{\tilde}{\widetilde}

\newcommand{\pt}{\{\bullet\}}
\renewcommand{\emptyset}{\varnothing}
\newcommand{\out}{\mathrm{out}}

\newtheorem{theorem}{Theorem}[section]
\newtheorem{proposition}[theorem]{Proposition}
\newtheorem{corollary}[theorem]{Corollary}

\theoremstyle{definition}
\newtheorem{definition}[theorem]{Definition}
\newtheorem{remark}[theorem]{Remark}
\newtheorem{example}[theorem]{Example}

\newtheorem*{ack}{Acknowledgements}
\numberwithin{equation}{section}

\begin{document}

\title{$2$-Segal sets and pseudomonoids in the bicategory of spans}
\author{Sophia E Marx}
\address{Department of Mathematics\\University of Massachusetts\\Amherst, MA}
\email{semarx@umass.edu}
\author{Rajan Amit Mehta}
\address{Department of Mathematical Sciences\\Smith College\\44 College Lane\\Northampton, MA 01063}
\email{rmehta@smith.edu}



\begin{abstract}
	In this survey article, we give an introduction to the notion of a $2$-Segal set and prove that $2$-Segal sets are equivalent to pseudomonoids in the bicategory of spans. The proof utilizes graphical techniques for $2$-Segal sets and spans that should be useful in more general settings. 

    There are procedures for obtaining an associative algebra from a $2$-Segal set (satisfying finiteness conditions). We describe these procedures and give several examples of algebras arising from $2$-Segal sets.

    Wherever possible, we avoid higher category theory so as to make the paper accessible to a wide audience.
\end{abstract}
\maketitle

\section{Introduction}\label{sec:introduction}

Given a category, one can form its \emph{nerve}, which is a simplicial set \cite{segal:classifying}. In the other direction, there are conditions, known as the \emph{Segal conditions}, that determine whether a given simplicial set is the nerve of a category. This correspondence leads to the more general notion of a \emph{Segal space} \cite{rezk:homotopy}, which is a simplicial space satisfying the Segal conditions up to weak homotopy equivalence. Intuitively, a Segal space can be viewed as a ``category up to coherent homotopy''.

The notion of \emph{$2$-Segal space} was introduced by Dyckerhoff and Kapranov \cite{Dyckerhoff-Kapranov:Higher}, and independently by G\'alvez-Carillo, Kock, and Tonks \cite{GKT1} (under the name \emph{decomposition space}), as a generalization of the notion of Segal space. Both groups of authors were motivated by the fact that one can, under certain conditions, obtain associative (co)algebras from $2$-Segal spaces; in \cite{Dyckerhoff-Kapranov:Higher}, the construction of Hall algebras was emphasized, and \cite{GKT1} was motivated by the study of incidence (co)algebras. 

Both groups also recognized that the $2$-Segal conditions provide higher coherence conditions for associativity. This idea was further developed by Stern \cite{Stern:span}, who considered the $\infty$-category of spans in $\cat$, where $\cat$ is an $\infty$-category with finite limits. He proved that there is an $\infty$-categorical equivalence between $2$-Segal objects in $\cat$ and coherent associative algebra objects in the $\infty$-category of spans in $\cat$.

In this review article, we consider Stern's result in the very special case where $\cat = \set$, the category of sets. In this case, the result (at the level of objects) is that there is an up-to-isomorphism correspondence between $2$-Segal sets and pseudomonoid objects in the bicategory of spans of sets. This is a statement that can be understood without any knowledge of $\infty$-category theory, and one might hope for an elementary proof that doesn't utilize $\infty$-category theory at all. This was partly done in \cite{CMS}, where one direction of the correspondence was described in detail, but the other direction was only roughly sketched. We provide here a full proof of the correspondence in both directions. In doing so, we make a couple of small improvements to the graphical approaches in \cite{CMS}.

The link between this result and the (co)algebra constructions of \cites{Dyckerhoff-Kapranov:Higher, GKT1} is given by the Baez-Dolan \emph{categorification} program \cites{BaezDolan:fromfinitesets, BaezGroupoid}, which we will briefly describe. ``Categorification'' is a general term used to describe situations where $n$-categorical data can be seen as arising from $(n+k)$-categorical data. A classical example of categorification is the fact that the Euler characteristic of a manifold, which is an element of the \emph{set} of integers, arises from its homology, which is an object in the \emph{category} of graded abelian groups. A further categorification comes from the recognition that homology comes from chain complexes, which are objects in the $2$-category of chain complexes, chain maps, and chain homotopies. This higher categorical data gives us a better understanding (for example, it explains why the Euler characteristic is invariant under homotopy equivalence), and it can also contain additional information (for example, homology reveals that $S^2$ is not homotopy equivalent to $S^4$, even though they both have Euler characteristic $2$). 

Functors out of (higher) span categories commonly arise in situations where structures can be both pulled back and pushed forward along maps, e.g.\ sheaves, and this provides a natural setting for categorification. An explicit example of this appears in \cite{BaezGroupoid}, where Baez, Hoffnung, and Walker defined a monoidal functor from the category of spans of groupoids (satisfying certain finiteness conditions) to the category $\vect_\kk$ of vector spaces over an arbitrary field $\kk$. 

In this paper, we are interested in the special case of sets (as opposed to groupoids), which simplifies the description of the functor considerably, and which we will now briefly describe; more details are given in Section \ref{sec:linearization}.

For any set $X$, let $\kk[X]$ be the vector space generated by elements of $X$. A map of sets $f: X \to Y$ naturally extends to a linear map $f_*: \kk[X] \to \kk[Y]$. Additionally, if the preimage $f^{-1}(y)$ is finite for all $y \in Y$, then we can define a map $f^*: \kk[Y] \to \kk[X]$, given by
\[ f^*(y) = \sum_{x \in f^{-1}(y)} x.\]

The category of spans of sets, which we denote $\Span_1$, is defined as follows. The objects are sets, the morphisms from $X$ to $Y$ are isomorphism classes of spans of the form
\begin{equation}\label{eqn:span}
    X \xleftarrow{f} A \xrightarrow{g} Y,
\end{equation} 
and composition of spans
\begin{equation}
    X \xleftarrow{} A \xrightarrow{} Y \xleftarrow{} B \xrightarrow{} Z
\end{equation} 
is given by taking the pullback of the diagram. 

Given a span \eqref{eqn:span} such that $f$ has finite preimages, we can obtain the map $g_* f^*: \kk[X] \to \kk[Y]$. This gives a well-defined functor from (a subcategory of) $\Span_1$ to $\vect_\kk$. Additionally, this functor is monoidal, taking the Cartesian product of sets and spans to the tensor product of vector spaces and linear maps.

The monoidal category $\Span_1$ can be viewed as the truncation of a monoidal bicategory $\Span$, where objects are sets, morphisms are spans, and $2$-morphisms are maps of spans; see Section \ref{sec:bicatspan} for more details. Thus, the bicategory $\Span$ provides a setting where categorification can occur. Specifically, given an associative algebra $\mathcal{A}$ over $\kk$, we can try to find a monoid object in $\Span_1$ whose image in $\vect_\kk$ is isomorphic to $\mathcal{A}$. We could then try to lift the monoid object to a coherent monoid structure in $\Span$. For bicategories, coherent monoid structures are commonly called \emph{pseudomonoids}, so a pseudomonoid in $\Span$ is precisely the structure that one should look for when categorifying algebras.

We can thus interpret Stern's correspondence as saying that a 2-Segal set is exactly the structure that provides the data of a categorified algebra, and as a result one can apply tools from the theory of simplicial sets to the study of algebraic structures. There is also a nice graphical calculus for $2$-Segal sets \cite{CMS} that helps to make the theory more accessible to readers from a variety of mathematical backgrounds and areas of interest.

We end the introduction with a few additional remarks:
\begin{itemize}
    \item Recent survey articles by Bergner \cite{bergner2024combinatorial} and Stern \cite{stern:perspectives} have some overlap with this one. We hope that the reader will find that this article complements theirs (and vice versa). In particular, Stern presents the correspondence between $2$-Segal sets and pseudomonoids in $\Span$ as an equivalence of categories, observing that morphisms of $2$-Segal sets correspond to oplax morphisms of pseudomonoids in $\Span$. On the other hand, our construction of a $2$-Segal set from a pseudomonoid in $\Span$ in Section \ref{sec:pseudomonoidto2Segal} is more explicit than Stern's, which makes use of the fact \cite{BOORS} that $2$-Segal sets are $3$-coskeletal.

    \item In \cite{CMS}, pseudomonoids in $\Span$ with additional structure were considered. There, it was shown that $2$-Segal paracyclic sets are categorifications of Frobenius algebras, and that $2$-Segal $\Gamma$-sets are categorifications of commutative algebras. 
    
    In \cite{MM:finset}, we show that, on a $2$-Segal set, a paracyclic and a $\Gamma$-structure will (without any additional compatibility conditions) fit together to form a \emph{cosymmetric set}, i.e.\ a functor $\Phi \to \set$, where $\Phi$ is the skeleton of the category of nonempty finite sets. It follows that $2$-Segal cosymmetric sets are categorifications of commutative Frobenius algebras, and thus can be viewed as categorified 2D topological quantum field theories.
    
    \item The middle step in the categorification process described above involves finding monoid objects in $\Span_1$. Such structures were considered in \cite{ContrerasMehtaSpan}, where it was shown that there is a correspondence between such monoid objects and simplicial sets that satisfy conditions that are similar to, but strictly weaker than, the $2$-Segal conditions. In \cite{CMS}*{Section 3.5}, it was observed that there are monoid objects in $\Span_1$ that do not lift to pseudomonoids in $\Span$, revealing that there are obstructions to categorification that may arise at this step.
    \item As noted above, Baez, Hoffnung, and Walker \cite{BaezGroupoid} considered the more general case of spans of groupoids. Some of their main examples (Hecke algebras and Hall algebras) appear in \cite{Dyckerhoff-Kapranov:Higher} as algebras that arise from $2$-Segal groupoids. The main examples in \cites{GKT1,GKT_comb} are also $2$-Segal groupoids. We plan to consider the case of $2$-Segal groupoids in more detail in future work.
    \item In \cite{Morton}, Morton constructed a $2$-functor from the bicategory of spans of (essentially finite) groupoids to that of Kapranov-Voevodsky $2$-vector spaces. Although he didn't prove that the functor respects monoidal structures, it seems likely that it does. If so, it should follow that one can obtain finitely semisimple $\kk$-linear monoidal categories from $2$-Segal sets. This seems to be closely related to the constructions in \cite{Dyckerhoff-Kapranov:Higher}*{Section 3.5}.
\end{itemize}

The structure of the paper is as follows: In Section \ref{sec:simplicial}, we review simplicial sets and Segal sets. In Section \ref{sec:2segal}, we define $2$-Segal sets in a way that closely parallels the definition of Segal sets. We also describe the graphical calculus for face and degeneracy maps in a $2$-Segal set. In Section \ref{sec:span}, we review the monoidal bicategory $\Span$ and the definition of pseudomonoid. In Section \ref{sec:correspondence}, we prove the correspondence between $2$-Segal sets and pseudomonoids in $\Span$. Finally, in Section \ref{sec:hallalgebra}, we describe the Hall algebra construction, as well as the closely related incidence (co)algebra construction, which produce associative (co)algebras from a $2$-Segal set. Throughout the paper we include numerous examples that may be of interest to the reader.

\begin{ack}
We thank Ivan Contreras and Walker Stern for helpful conversations about topics related to the paper. We also thank the organizers of the ECOGyT conference in Bogot\'a in 2024 for the inspiration to write this paper. Finally, we thank the anonymous referee for many helpful suggestions and corrections.
\end{ack}

\section{Simplicial sets, categories, and the Segal conditions}
\label{sec:simplicial}

\subsection{Simplicial sets}
We begin with a brief exposition of simplicial sets, since notation and terminology for these tend to vary even among standard sources (see e.g. \cite{goerss-jardine} for a more thorough treatment).

\begin{definition}
    The \emph{simplex category} $\Delta$ is the category whose objects are the ordered sets $[n] = \{0 < 1 < \dots < n\}$ for $n = 0,1,2, \dots$, and whose morphisms are weakly order-preserving set maps $f: [m] \rightarrow [n]$.
\end{definition}

The morphisms of $\Delta$ are generated by the injections $\delta_i^n: [n-1] \to [n]$ for $0 \leq i \leq n$, 
\[ \delta_i^n(k) = \begin{cases}
    k, & k < i,\\
    k+1 & k \geq i,
\end{cases}
\]
called the \emph{coface maps},
and the surjections $\sigma_i^n: [n+1] \to [n]$ for $0 \leq i \leq n$,
\[ \sigma_i^n(k) = \begin{cases}
    k, & k \leq i, \\
    k-1, & k > i,
\end{cases}
\]
called the \emph{codegeneracy maps}. They satisfy the relations
\begin{align}
    \delta_j^{n} \delta_i^{n-1} &= \delta_i^{n} \delta_{j-1}^{n-1},  \hspace{28px} i < j, \label{eqn:coface1}\\
    \sigma_j^{n} \delta_i^{n-1} &= \begin{cases}
        \delta_i^{n+1} \sigma_{j-1}^n, & i < j,\\
        \id, &  i=j, j+1,\\
        \delta_{i-1}^{n+1} \sigma_j^n, &  i > j+1,
    \end{cases}\\
    \sigma_j^{n} \sigma_i^{n+1} &= \sigma_i^{n} \sigma_{j+1}^{n+1}, \hspace{26px} i \leq j. \label{eqn:coface3}
\end{align}

\begin{definition}
    
A \emph{simplicial set} is a functor $X_\bullet: \Delta^\op \to \set$. 

\end{definition}

Using the above description of $\Delta$, we can describe a simplicial set via the following data:
\begin{itemize}
    \item sets $X_n$ for $n=0,1,2,\dots$,
    \item \emph{face maps} $d_i^n: X_n \to X_{n-1}$ for $0 \leq i \leq n$, and
    \item \emph{degeneracy maps} $s_i^n: X_n \to X_{n+1}$ for $0 \leq i \leq n$.
\end{itemize}
This data must satisfy the \emph{simplicial relations} dual to \eqref{eqn:coface1}--\eqref{eqn:coface3}:
\begin{align}
    d_i^{n-1}d_j^n &= d_{j-1}^{n-1} d_i^n, \hspace{28px} i < j, \label{eqn:face1}\\
    d_i^{n-1} s_j^n &= \begin{cases}
    s_{j-1}^n d_i^{n+1}, & i < j,\\
    \id, & i=j,j+1,\\
    s_j^n d_{i-1}^{n+1}, & i > j+1,
    \end{cases}\\
    s_i^{n+1} s_j^n &= s_{j+1}^{n+1} s_i^n, \hspace{28px} i \leq j.\label{eqn:face3}
\end{align}

Intuitively, we can think of an element of $X_n$ as being an $n$-simplex with vertices labeled $\{0, \dots, n\}$. Then the face map $d_i^n$ can be thought of as the map that picks out the face opposite to the vertex $i$, and the degeneracy map $s_i^n$ can be thought of as the map that views an $n$-simplex as an $(n+1)$-simplex where vertices $i$ and $i+1$ coincide.

\begin{definition}
    Given a simplicial set $X_\bullet$, we define the following special maps:
\begin{itemize}
    \item the \emph{vertex maps} $v_i^n : X_n \to X_0$, $0 \leq i \leq n$, are those induced by the map $[0] \to [n]$, $0 \mapsto i$,
    \item the \emph{vertebra maps} $e_i^n: X_n \to X_1$, $1 \leq i \leq n$, are induced by the map $[1] \to [n]$, $0 \mapsto i-1$, $1 \mapsto i$,
    \item the \emph{long edge map} $e_\out^n: X_n \to X_1$, is induced by the map $[1] \to [n]$, $0 \mapsto 0$, $1 \mapsto n$.
\end{itemize}

\end{definition}

All of these maps can be expressed (not uniquely) in terms of face and degeneracy maps. In particular,
\begin{align*}
    v_0^0 &= \id_{X_0}, & v_0^1 &= d_1^1, & v_1^1 &= d_0^1, \\
    e_1^1 &= \id_{X_1}, & e_1^2 &= d_2^2, & e_2^2 &= d_0^2,\\
    e_\out^0 &= s_0^0, & e_\out^1 &= \id_{X_1}, & e_\out^2 &= d_1^2,
\end{align*}
and, for example,
\begin{align*}
    v_1^3 &= d_0^1 d_2^2 d_3^3 = d_1^1 d_1^2 d_0^3,\\
    e_2^4 &= d_0^2 d_3^3 d_4^4 = d_2^2 d_2^3 d_0^4, \\
    e_\out^3 &= d_1^2 d_1^3 = d_1^2 d_2^3.
\end{align*}
Intuitively, $v_i^n$ can be thought of as the map that picks out the $i$th vertex of an $n$-simplex. Similarly, $e_i^n$ picks out the edge connecting vertices $i-1$ and $i$, and $e_\out^n$ picks out the edge connecting vertices $0$ and $n$.

We will often suppress the upper index in face, degeneracy, vertex, vertebra, and long edge maps when it is clear from context.

\subsection{The nerve of a category}\label{sec:nervecat}
Let $\cat$ be a small category. Then one can form a simplicial set called the \emph{nerve} of $\cat$, defined as follows.

The first step is to recognize that $\Delta$ can be embedded in $\Cat$, the category of (small) categories, by viewing each $[n]$ as a poset category, i.e.\ the category with objects $0,\dots,n$, and with a unique morphism from $i$ to $j$ whenever $i \leq j$. 

\begin{definition}
    The \emph{nerve} of a category $\cat$, denoted $N\cat_\bullet$, is given by 
    \[N\cat_n = \Hom_{\Cat} ([n], \cat).\]
\end{definition}
This definition is expedient, but it would be nice to have a more explicit description. To do this, we observe that an element of $N\cat_n$ consists of the following data:
\begin{itemize}
    \item objects $x_0, \dots, x_n$ of $\cat$,
    \item morphisms $f_{ij} \in \Hom_\cat(x_i,x_j)$ for $0 \leq i \leq j \leq n$,
\end{itemize}
such that $f_{ik} = f_{jk} \circ f_{ij}$ whenever $i \leq j \leq k$. From this property, one can see that such data is uniquely determined by the morphisms $f_{01}, f_{12}, \dots, f_{(n-1)n}$. Conversely, any $n$-tuple $(f_{01}, f_{12}, \dots, f_{(n-1)n})$ of morphisms in $\cat$ that are composable, in the sense that the source of $f_{i(i+1)}$ equals the target of $f_{(i-1)i}$, uniquely extends to an element of $N\cat_n$.

The conclusion is that we can identify $N\cat_n$ with the set of composable $n$-tuples $(f_1, f_2, \dots, f_n)$, where we write $f_i$ as shorthand for $f_{(i-1)i}$. Under this identification, the face maps are given by
\begin{align*}
    d_0(f_1, \dots, f_n) &= (f_2, \dots, f_n),\\
    d_i(f_1, \dots, f_n) &= (f_1, \dots, f_{i+1}\circ f_i, \dots, f_n) \mbox{ for } 0<i<n, \\
    d_n(f_1, \dots, f_n) &= (f_1, \dots, f_{n-1}), 
\end{align*}
and the degeneracy maps are given by insertion of identity morphisms.

\subsection{Segal sets}\label{sec:segal}
The \emph{Segal conditions} are conditions that characterize the simplicial sets that arise as nerves of categories. They are due to Grothendieck \cite{grothendieck:FGA3}*{Section 4}, and the naming is due to their (very brief!) appearance in \cite{segal:classifying}. In this subsection we will introduce the Segal conditions.

\begin{definition}\label{def:piqi}
    For $0 < i < n$, the inclusion maps
\begin{align*}
    P_i^n&: [n-i] \to [n], \\
    Q_i^n&: [i] \to [n],
\end{align*}
are defined by $P_i^n(k) = k+i$ and $Q_i^n(k) = k$. 
\end{definition}
The maps in Definition \ref{def:piqi} can be visualized by considering the interval $[0,n]$ as representing the ordered set $[n]$; see Figure \ref{fig:interval}. Subdividing the interval at the point $i$ creates the subintervals $[0,i]$ and $[i,n]$, which are the images of the maps $Q_i^n$ and $P_i^n$, respectively. The fact that the two subintervals share the common point $i$ reflects the fact that there is a commutative diagram in $\Delta$:
\begin{equation}\label{diag:intervalsubdivide}
\begin{tikzcd}
{[0]} \arrow[r, "0"] \arrow[d, "i"']       & {[n-i]} \arrow[d, "P_{i}^n"] \\
{[i]} \arrow[r, "Q_{i}^n"] & {[n]}                      
\end{tikzcd}
\end{equation}

\begin{figure}[t]
\begin{center}
\begin{tikzpicture}[scale=0.4]
    \begin{scope} [shift={(0,0)}]
        \draw (0,0) node[vertex] (0) {} node[anchor=north] {0}
        -- (1,0) node[vertex] (1) {} node[anchor=north] {1}
        -- (2,0) node[vertex] (2) {} node[anchor=north] {2}
        -- (2.5,0) node[vertex] (3) {} node[anchor=north] {};

        \node at (3.4,0) {$\cdots$};
        
       \draw (4.2,0) -- (4.7,0) node[vertex] (i) {} node[anchor=north] {$i$} -- (5.2,0);
       
        \node at (6.1,0) {$\cdots$};
        
       \draw (7,0) -- (7.5,0) node[vertex] (n-1) {} node[anchor=north] {}
       -- (8.5,0) node[vertex] (n) {} 
       node[anchor=north] {$n$};
       
        \draw[fill] (0) circle (2pt);
        \draw[fill] (1) circle (2pt);
        \draw[fill] (2) circle (2pt);
        \draw[fill] (i) circle (2pt);
        \draw[fill] (n-1) circle (2pt);
        \draw[fill] (n) circle (2pt);

    \end{scope}
    \begin{scope} [shift={(0,-3)}]
        \draw (0,0) node[vertex] (0) {} node[anchor=north] {0}
        -- (1,0) node[vertex] (1) {} node[anchor=north] {1}
        -- (2,0) node[vertex] (2) {} node[anchor=north] {2}
        -- (2.5,0) node[vertex] (3) {} node[anchor=north] {};

        \node at (3.4,0) {$\cdots$};
        
       \draw (4.2,0) -- (4.7,0) node[vertex] (i) {}; 
       \draw  (5.2,0) node[vertex] (i+1) {} 
       -- (5.7,0) ;
       \node at (4.95,0) [anchor=north] {$i$};

        \node at (6.6,0) {$\cdots$};
        
       \draw (7.5,0) -- (8,0) node[vertex] (n-1) {} node[anchor=north] {}
       -- (9,0) node[vertex] (n) {} 
       node[anchor=north] {$n$};
       
        \draw[fill] (0) circle (2pt);
        \draw[fill] (1) circle (2pt);
        \draw[fill] (2) circle (2pt);
        \draw[fill] (i) circle (2pt);
        \draw[fill] (i+1) circle (2pt);
        \draw[fill] (n-1) circle (2pt);
        \draw[fill] (n) circle (2pt);

    \end{scope}
\end{tikzpicture}
\end{center}
\caption{The interval $[0,n]$, representing the ordered set $[n]$. Below it, the interval is subdivided at the point $i$. The resulting subintervals, $[0,i]$ and $[i,n]$, are the images of the maps $Q_i^n$ and $P_i^n$, respectively.}
\label{fig:interval}
\end{figure}

Given a simplicial set $X_\bullet$, we obtain from diagram \ref{diag:intervalsubdivide} a commutative diagram of sets
\begin{equation}\label{diag:intervalsubdivide2}
\begin{tikzcd}
X_n \arrow[r, "\hat{Q}_{i}^n"] \arrow[d, "\hat{P}_{i}^n"'] & X_{i} \arrow[d, "v_i"] \\
X_{n-i} \arrow[r, "v_0"]                            & X_0
\end{tikzcd}
\end{equation}
and thus a map 
\begin{equation}\label{eqn:intervalsubdivide0}
X_n \to X_i \times_{X_0} X_{n-i}.    
\end{equation}

More generally, for $0 < i_1 < \cdots < i_\ell < n$, we can produce a subdivision of the interval $[0,n]$ into the subintervals $[0, i_1], [i_1,i_2], \dots, [i_{\ell-1},i_\ell], [i_\ell, n]$. Such a subdivision induces a map
\begin{equation}\label{eqn:intervalsubdivide}
    X_n \to X_{k_0} \times_{X_0} X_{k_1} \times_{X_0} \cdots \times_{X_0} X_{k_\ell},
\end{equation}
where $k_0 = i_1$, $k_p = i_{p+1} - i_p$ for $0 < p < \ell$, and $k_\ell = n - i_\ell$. We note that the map \eqref{eqn:intervalsubdivide} can be constructed either in a single step, as a generalization of \eqref{eqn:intervalsubdivide0}, or by iterative application of \eqref{eqn:intervalsubdivide0}.

In particular, we can consider the finest subdivision $0 <1 <2 <\cdots < n$, from which we obtain the map
\begin{equation}\label{eqn:spine}
\begin{split}
T_n: X_n &\to \underbrace{X_1 \times_{X_0} X_1 \times_{X_0} \cdots \times_{X_0} X_1}_\text{n factors}, \\
\omega &\mapsto (e_1 \omega,\dots, e_n \omega).
\end{split}
\end{equation}
The map $T_n$ in \eqref{eqn:spine} is often referred to as the \emph{spine map} or the \emph{Segal map}. This explains the terminology ``vertebra maps'', since they are the components of the spine map.

\begin{proposition}\label{prop:segal}
    Let $X_\bullet$ be a simplicial set. The following are equivalent:
    \begin{enumerate}
        \item The map \eqref{eqn:intervalsubdivide0} is an isomorphism for all $0 < i < n$.
        \item The map \eqref{eqn:intervalsubdivide} is an isomorphism for all $0 < i_1 < \cdots < i_\ell < n$.
        \item The spine map \eqref{eqn:spine} is an isomorphism for all $n$.
    \end{enumerate}
\end{proposition}

\begin{proof}
    The implications $2 \implies 1$ and $2 \implies 3$ are immediate, since the maps \eqref{eqn:intervalsubdivide0} and \eqref{eqn:spine} are special cases of \eqref{eqn:intervalsubdivide}. The implication $1 \implies 2$ follows from the fact that the maps in \eqref{eqn:intervalsubdivide} can be obtained by repeated use of \eqref{eqn:intervalsubdivide0}. Finally, the implication $3 \implies 2$ follows from the commutative diagram
    \begin{equation*}
    \begin{tikzcd}
X_n \arrow[r] \arrow[rd, swap, "T_n"] & X_{k_0} \times_{X_0} \cdots \times_{X_0} X_{k_\ell} \arrow[d, "T_{k_0} \times \cdots \times T_{k_\ell}"] \\
 & X_1 \times_{X_0} \cdots \times_{X_0} X_1
\end{tikzcd}
\end{equation*}
which shows that, if all spine maps are isomorphisms, then the map \eqref{eqn:intervalsubdivide} is an isomorphism.
\end{proof}

\begin{definition}
A simplicial set $X_{\bullet}$ is \emph{Segal} (or \emph{$1$-Segal}) if any of the equivalent conditions in Proposition \ref{prop:segal} are satisfied.
\end{definition}

\subsection{The Segal functors}

In this subsection, we will describe yet another perspective on Segal sets, which will be useful in proving the following key fact:
\begin{theorem}[\cite{grothendieck:FGA3}*{Proposition 4.1}]\label{thm:1segal}
A simplicial set is Segal if and only if it is isomorphic to the nerve of a category.
\end{theorem}

Let $X_\bullet$ be a simplicial set. As we saw in Section \ref{sec:segal}, there is a map \eqref{eqn:intervalsubdivide} associated to each subdivision of the interval $[0,n]$. These maps can be composed in certain ways, since a subinterval can itself be subdivided. Taken together, the subdivision maps and their compositions form a nice structure that we will now describe.

For each $n \geq 1$, consider the poset of subdivisions $0 < i_1 < \cdots < i_\ell < n$, where finer subdivisions are ``greater than'' coarser subdivisions. The Hasse diagram of this poset has the structure of the $(n-1)$-cube, and in particular there is an initial object (the unsubdivided interval) and a terminal object (the fully subdivided interval). The first nontrivial cases, $n=2$ and $n=3$, are shown in Figure \ref{fig:bin2}.
\begin{figure}[t]
\begin{center}
\begin{tikzpicture}[scale=0.4]
    \begin{scope} [shift={(0,-.25)}]
        \draw (0,0) node[vertex] (0) {} node[anchor=north] {0}
        -- (1,0) node[vertex] (1) {} node[anchor=north] {1}
        -- (2,0) node[vertex] (2) {} node[anchor=north] {2};
        \draw[fill] (0) circle (2pt);
        \draw[fill] (1) circle (2pt);
        \draw[fill] (2) circle (2pt);
    \end{scope}
    \begin{scope} [shift={(-.2,-3.75)}]
        \draw (0,0) node[vertex] (0) {} node[anchor=north] {0}
        -- (1,0);
        \draw (1.4,0) 
        -- (2.4,0) node[vertex] (2) {} node[anchor=north] {2};
    \node at (1.2,0) [anchor=north] {1};
        \draw[fill] (0) circle (2pt);
        \draw[fill] (2) circle (2pt);
        \draw[fill] (1,0) circle (2pt);
        \draw[fill] (1.4,0) circle (2pt);
    \end{scope}
    \node[rotate=270] at (1,-2.25) {$\Longrightarrow$};    
    \begin{scope} [shift={(15,0)}]
            \begin{scope} [shift={(0,0)}]
        \draw (0,0) node[vertex] (0) {} node[anchor=north] {0}
        -- (1,0) node[vertex] (1) {} node[anchor=north] {1}
        -- (2,0) node[vertex] (2) {} node[anchor=north] {2}
        -- (3,0) node[vertex] (3) {} node[anchor=north] {3};
        \draw[fill] (0) circle (2pt);
        \draw[fill] (1) circle (2pt);
        \draw[fill] (2) circle (2pt);
        \draw[fill] (3) circle (2pt);
    \end{scope}
    \begin{scope} [shift={(-7,-2)}]
        \draw (0,0) node[vertex] (0) {} node[anchor=north] {0}
        -- (1,0);
        \draw (1.4,0) 
        -- (2.4,0) node[vertex] (2) {} node[anchor=north] {2}        
        -- (3.4,0) node[vertex] (3) {} node[anchor=north] {3};
    \node at (1.2,0) [anchor=north] {1};
        \draw[fill] (0) circle (2pt);
        \draw[fill] (1,0) circle (2pt);
        \draw[fill] (1.4,0) circle (2pt);
        \draw[fill] (2) circle (2pt);
        \draw[fill] (3) circle (2pt);
    \end{scope}    
    \begin{scope} [shift={(6.6,-2)}]
        \draw (0,0) node[vertex] (0) {} node[anchor=north] {0}
        -- (1,0) node[vertex] (1) {} node[anchor=north] {1}
        -- (2,0);
        \draw (2.4,0)
        -- (3.4,0) node[vertex] (3) {} node[anchor=north] {3};
        \node at (2.2,0) [anchor=north] {2};
        \draw[fill] (0) circle (2pt);
        \draw[fill] (1) circle (2pt);
        \draw[fill] (2,0) circle (2pt);
        \draw[fill] (2.4,0) circle (2pt);
        \draw[fill] (3) circle (2pt);
    \end{scope}
    \begin{scope} [shift={(-.4,-4)}]
        \draw (0,0) node[vertex] (0) {} node[anchor=north] {0}
        -- (1,0);
        \draw (1.4,0)
        -- (2.4,0);
        \draw (2.8,0)
        -- (3.8,0) node[vertex] (3) {} node[anchor=north] {3};
        \node at (1.2,0) [anchor=north] {1};
        \node at (2.6,0) [anchor=north] {2};
        \draw[fill] (0) circle (2pt);
        \draw[fill] (1,0) circle (2pt);
        \draw[fill] (1.4,0) circle (2pt);        
        \draw[fill] (2.4,0) circle (2pt);
        \draw[fill] (2.8,0) circle (2pt);
        \draw[fill] (3) circle (2pt);
    \end{scope}  
    \node[rotate=330] at (5,-1) {$\Longrightarrow$};
    \node[rotate=210] at (-2,-1) {$\Longrightarrow$};
    \node[rotate=210] at (5,-3) {$\Longrightarrow$};
    \node[rotate=330] at (-2,-3) {$\Longrightarrow$};
    \end{scope}
\end{tikzpicture}
\end{center}
\caption{The Hasse diagrams of the posets of subdivisions of the interval $[0,n]$ for $n=2$ (on the left) and $n=3$ (on the right).}
\label{fig:bin2}
\end{figure}

\begin{figure}[t]
\begin{center}
\begin{tikzcd}
X_2 \arrow[d, "T_2"] \\ X_1 \times_{X_0} X_1
\end{tikzcd}
\hspace{.5cm}
\begin{tikzcd}[column sep=small]
& X_3 \arrow[ld, "{(e_1,d_0)}"'] \arrow[rd, "{(d_3,e_3)}"] \arrow[dd,"T_3"] & \\
X_1 \times_{X_0} X_2 \arrow[rd, "\id \times T_2"'] & & X_2 \times_{X_0} X_1 \arrow[ld, "T_2 \times \id"] \\
& X_1 \times_{X_0} X_1 \times_{X_0} X_1 &
\end{tikzcd}
\end{center}
    \caption{Images of the Segal functors for $n=2$ and $n=3$.}
    \label{fig:segalfunctor}
\end{figure}

Collectively, the maps \eqref{eqn:intervalsubdivide} can be viewed as giving, for each $n \geq 1$, a functor from the poset category of subdivisions of $[0,n]$ to the category of sets. The commutative diagrams showing the images of this functor for $n=2$ and $n=3$ are given in Figure \ref{fig:segalfunctor}. We call these functors the \emph{Segal functors} associated to $X_\bullet$. In particular, the Segal functors take the initial-to-terminal morphisms in the subdivision posets to the spine maps \eqref{eqn:spine}.

As an immediate consequence of the definitions, we have
\begin{proposition}\label{prop:1segalfunctor}
    A simplicial set $X_\bullet$ is Segal if and only if, for all $n \geq 1$, the image of the associated Segal functor consists of isomorphisms.
\end{proposition}

We are now ready to prove Theorem \ref{thm:1segal}.

\begin{proof}[Proof of Theorem \ref{thm:1segal}]
The ``if'' direction follows from the observation that the identification of $N\cat_n$ with the set of composable $n$-tuples in Section \ref{sec:nervecat} is exactly the spine map.

For the ``only if'' direction, suppose that $X_\bullet$ is Segal. We construct a category $\cat$ as follows.
\begin{itemize}
    \item The objects of $\cat$ are the elements of $X_0$.
    \item Given two objects $x,y \in X_0$, a morphism from $x$ to $y$ is an element $\alpha \in X_1$ such that $d_1 \alpha = x$ and $d_0 \alpha  = y$. Thus, $X_1$ is the set of all morphisms in $\cat$.
    \item Given $\beta \in \Hom_\cat(x,y)$ and $\alpha \in \Hom_\cat(y,z)$, there is a unique $\omega \in X_2$ such that $T_2 \omega  = (\beta, \alpha)$. Then the composition in $\cat$ of $\alpha$ and $\beta$ is defined as $\alpha \circ \beta = d_1 \omega$.
\end{itemize}

To see that the composition in $\cat$ is associative, consider the following diagram, which commutes as a result of the simplicial identities. The arrows in bold represent maps in the image of the Segal functors (see Figure \ref{fig:segalfunctor}) and are therefore isomorphisms.
\begin{equation*}
\begin{tikzcd}[column sep=tiny]
& & X_1 \times_{X_0} X_1 \times_{X_0} X_1 &  & \\
& X_2 \times_{X_0} X_1 \arrow[ru, LA, "T_2 \times \id"] \arrow[ld, "d_1 \times \id"'] & X_3 \arrow[u, LA, "T_3"] \arrow[r, LA, "{(e_1,d_0)}"] \arrow[l, LA, "{(d_3,e_3)}"'] \arrow[ld, "d_1"] \arrow[rd, "d_2"'] & X_1 \times_{X_0} X_2 \arrow[lu, LA, "\id \times T_2"'] \arrow[rd, "\id \times d_1"] & \\ 
X_1 \times_{X_0} X_1 & X_2 \arrow[l, LA, "T_2"'] \arrow[rd, "d_1"]  &  & X_2 \arrow[r, LA, "T_2"] \arrow[ld, "d_1"'] & X_1 \times_{X_0} X_1 \\
& & X_1& &
\end{tikzcd}
\end{equation*}
Starting with $(\gamma,\beta,\alpha) \in X_1 \times_{X_0} X_1 \times_{X_0} X_1$, proceeding along the left side of the diagram to $X_1$ produces $\alpha \circ (\beta \circ \gamma) \in X_1$, and proceeding along the right edge produces $(\alpha \circ \beta) \circ \gamma$, so associativity follows from commutativity of the diagram.

To see that identity morphisms exist in $\cat$, consider the following commutative diagram.
\begin{equation*}
\begin{tikzcd}
X_1 \arrow[r, "s_0"] \arrow[rd, "\id"] & X_2 \arrow[r, LA, "T_2"] \arrow[d, "d_1"] & X_1 \times_{X_0} X_1 \\ & X_1 &
\end{tikzcd}
\end{equation*}
Given $\alpha \in \Hom_\cat(x,y)$, we have that
\[ T_2(s_0 \alpha) = (d_2 s_0 \alpha, d_0 s_0 \alpha) = (s_0 d_1 \alpha,\alpha) = (s_0 x, \alpha),\]
and it follows that
\[ \alpha \circ s_0 x = d_1 s_0 \alpha = \alpha.\]
A similar calculation using $s_1 \alpha$ shows that $s_0 y \circ \alpha = \alpha$. This establishes that $\cat$ is a category, but we still need to prove that $X_\bullet$ is isomorphic to $N\cat_\bullet$.

Using the identification of $N\cat_n$ with the set of composable $n$-tuples of morphisms (see Section \ref{sec:nervecat}), we immediately have isomorphisms of sets
\[ X_n \xrightarrow{T_n} \underbrace{X_1 \times_{X_0} X_1 \times_{X_0} \cdots \times_{X_0} X_1}_\text{n factors} \cong N\cat_n\]
for $n \geq 1$, with $X_0 = N\cat_0$ by construction. It remains to show that these isomorphisms are compatible with face and degeneracy maps. Because the spine maps can be expressed in terms of face maps, this can be shown in a straightforward way using the simplicial relations.
\end{proof}

\section{2-Segal sets}
\label{sec:2segal}

Geometrically, a Segal set is $1$-dimensional, in the sense that an $n$-simplex is uniquely determined by its spine, consisting of vertebrae that are $1$-simplices. A $2$-Segal set can be viewed as a $2$-dimensional analogue, where an $n$-simplex is uniquely determined by a collection of $2$-simplices.

There are several equivalent ways to define $2$-Segal sets. We will give three such definitions, analogous to the three definitions of $1$-Segal set given in Proposition \ref{prop:segal}. We will then describe a graphical calculus for $2$-Segal sets that appeared in \cite{CMS}.

\subsection{Definition of 2-Segal set}\label{sec:2segaldefinition}

\begin{definition} \label{def:PQ}
    
For $0 \leq i \leq j \leq n$, the maps 
\begin{align*}
  P_{ij}^n&: [j-i] \to [n],\\
  Q_{ij}^n&: [n+i-j+1] \to [n]
\end{align*}
are defined by
\begin{align*}
P_{ij}^n(k) &= k+i,\\
    Q_{ij}^n(k) &= \begin{cases}
        k, & 0 \leq k \leq i, \\
        k+j-i-1, & i < k \leq n+i-j+1.
    \end{cases}
\end{align*}

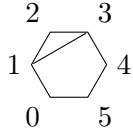
\begin{figure}[t]
    \begin{center}
    \begin{tikzpicture}[scale=0.5]
            \begin{scope}
        \draw (240:1) node[vertex] (0) {} node[anchor=north east] {0}
        -- (300:1) node[vertex] (5) {} node[anchor=north west] {5}
        -- (0:1) node[vertex] (4) {} node[anchor=west] {4}
        -- (60:1) node[vertex] (3) {} node[anchor=south west] {3}
        -- (120:1) node[vertex] (2) {} node[anchor=south east] {2}
        -- (180:1) node[vertex] (1) {} node[anchor=east] {1}        
        -- cycle;
        \draw (1) -- (3);
    \end{scope}
    \end{tikzpicture}
    \end{center}
    \caption{Visualization of the maps $P_{13}^5: [2] \to [5]$, given by $\{0,1,2\} \mapsto \{1,2,3\}$, and $Q_{13}^5: [4] \to [5]$, given by $\{0,1,2,3,4\} \mapsto \{0,1,3,4,5\}$.}
    \label{fig:edgedivision}
\end{figure}

\end{definition}

The maps in Definition \ref{def:PQ} can be visualized by considering a regular $(n+1)$-gon with the vertices labeled in order from $0$ to $n$. The line segment connecting the vertices $i$ and $j$ subdivides the $(n+1)$-gon into a $(j-i+1)$-gon and an $(n+i-j+2)$-gon; see Figure \ref{fig:edgedivision}. The fact that the two polygons in the subdivision share a common edge reflects the fact that there is a commutative diagram in $\Delta$:
\begin{equation}\label{diag:subdivide}
\begin{tikzcd}
{[1]} \arrow[r, "{\{0,j-i\}}"] \arrow[d, "{\{i,i+1\}}"']       & {[j-i]} \arrow[d, "P_{ij}^n"] \\
{[n+i-j+1]} \arrow[r, "Q_{ij}^n"] & {[n]}                      
\end{tikzcd}
\end{equation}
We remark that the maps $P_{ij}^n$ and $Q_{ij}^n$, as well as the commutative diagram in equation \eqref{diag:subdivide}, are still well-defined in cases where the interpretation as polygonal subdivision degenerates, i.e.\ when $j-i = 0,1$, or $n$.

Given a simplicial set $X_\bullet$, the diagram \eqref{diag:subdivide} induces a commutative diagram of sets
\begin{equation}\label{diag:subdivide2}
\begin{tikzcd}
X_n \arrow[r, "\hat{Q}_{ij}^n"] \arrow[d, "\hat{P}_{ij}^n"'] & X_{n+i-j+1} \arrow[d, "e_{i+1}"] \\
X_{j-i} \arrow[r, "e_{\out}"]                            & X_1
\end{tikzcd}
\end{equation}
and thus a map
\begin{equation}\label{eqn:subdivide0}
    X_n \to X_{j-i} \times_{X_1} X_{n+i-j+1}.
\end{equation}

As mentioned previously, the diagram \eqref{diag:subdivide2} can be seen (except in some degenerate cases) as arising from a subdivision of the $(n+1)$-gon into two polygons. More generally, given a subdivision of the $(n+1)$-gon into $(k_1+1), \dots, (k_\ell+1)$-gons, the maps in \eqref{diag:subdivide2} can be repeatedly used to produce a map
\begin{equation}\label{eqn:subdivide3}
    X_n \to X_{k_1} \times_{X_1} \cdots \times_{X_1} X_{k_\ell}.
\end{equation}
In particular, we can consider a triangulation of the $(n+1)$-gon, from which we obtain a map
\begin{equation}\label{eqn:triangulation}
X_n \to \underbrace{X_2 \times_{X_1} X_2 \times_{X_1} \dots \times_{X_1} X_2}_\text{n-1 factors}. 
\end{equation}

\begin{proposition}\label{prop:2segalequivalent}
Let $X_\bullet$ be a simplicial set. The following are equivalent:
\begin{enumerate}
    \item The map \eqref{eqn:subdivide0} is an isomorphism for all $0 \leq i < j \leq n$.
    \item The map \eqref{eqn:subdivide3} is an isomorphism for all subdivided $(n+1)$-gons.
    \item The map \eqref{eqn:triangulation} is an isomorphism for all triangulated $(n+1)$-gons.
\end{enumerate}
\end{proposition}

The proof of Proposition \ref{prop:2segalequivalent} is similar to that of Proposition \ref{prop:segal}; see \cite{Dyckerhoff-Kapranov:Higher}*{Proposition 2.3.2}, and see Remark \ref{remark:unitality} below for a discussion of the degenerate cases of \eqref{eqn:subdivide0}.

\begin{definition}
We say that $X_\bullet$ is \emph{$2$-Segal} if any of the equivalent conditions in Proposition \ref{prop:2segalequivalent} are satisfied.
\end{definition}

\begin{remark}\label{remark:unitality}
In the case where $i=j$, \eqref{diag:subdivide2} becomes
\begin{equation}\label{diag:unital}
\begin{tikzcd}
X_n \arrow[r, "s_i"] \arrow[d, "v_i"'] & X_{n+1} \arrow[d, "e_{i+1}"] \\
X_{0} \arrow[r, "s_0"]                            & X_1
\end{tikzcd}
\end{equation}
In \cite{Dyckerhoff-Kapranov:Higher}, the term \emph{unital} was used to describe $2$-Segal spaces for which \eqref{diag:unital} is a pullback for all $n$. However, it was shown in \cite{FGKW:unital} that every $2$-Segal space is unital; in other words, \eqref{diag:unital} is automatically a pullback as a consequence of the conditions in Proposition \ref{prop:2segalequivalent}. This result might seem counterintuitive, but we note that the proof in \cite{FGKW:unital} relies on properties of the degeneracy maps in $X_\bullet$.

In the other degenerate cases, where $i+1=j$ or $i=0,j=n$, \eqref{diag:subdivide2} becomes
\begin{equation}
\begin{tikzcd}
X_n \arrow[r, "\id"] \arrow[d, "e_{i+1}"'] & X_n \arrow[d, "e_{i+1}"] \\
X_1 \arrow[r, "\id"]                            & X_1
\end{tikzcd}
\hspace{1cm}\mbox{or}\hspace{1cm}
\begin{tikzcd}
X_n \arrow[r, "e_\out"] \arrow[d, "\id"'] & X_1 \arrow[d, "\id"] \\
X_n \arrow[r, "e_\out"]                            & X_1
\end{tikzcd}
\end{equation}
which are automatically pullbacks.

\end{remark}

\subsection{The nerve of a partial category}
A fairly broad class of examples of $2$-Segal sets arise from a structure that we call a \emph{partial category}. Intuitively, a partial category is like a category, except composition of morphisms is not always defined, even when their sources and targets match up appropriately. 

This concept has appeared before under the name \emph{precategory} \cite{mateus-sernadas}, but we avoid this term because it is also used in other contexts. Partial categories can be seen as special cases of \emph{paracategories} \cites{hermida-mateus1, hermida-mateus2}, and they are a natural generalization of the notion of \emph{partial monoid} \cite{Segal:Conf}. It was shown in \cite{BOORS} that nerves of partial monoids are $2$-Segal, and we show here that this result holds more generally for partial categories\footnote{Partial categories can also be viewed as a special case of the $\mu$-categories of \cite{Dyckerhoff-Kapranov:Higher}*{Definition 3.3.3}, giving an alternative approach to proving that their nerves are $2$-Segal.}.

\begin{definition}\label{def:partialcat}
A (small) \emph{partial category} $\cat$ consists of
\begin{itemize}
    \item a set of objects $\Ob(\cat)$,
    \item for each $x,y \in \Ob(\cat)$, a set $\cat_{x,y}$ of morphisms from $x$ to $y$,
    \item for each $x,y,z \in \Ob(\cat)$, a set $\cat_{x,y,z} \subseteq \cat_{x,y} \times \cat_{y,z}$ with a composition map $\circ: \cat_{x,y,z} \to \cat_{x,z}$, $(f,g) \mapsto g \circ f$,
\end{itemize}
such that
\begin{enumerate}
    \item (associativity) for all $f \in \cat_{y,z}$, $g \in \cat_{x,y}$, $h \in \cat_{w,x}$, 
    \[(f \circ g) \circ h = f \circ (g \circ h).\] 
    \item (unitality) for all $x \in \Ob(\cat)$, there exists $\id_x \in \cat_{x,x}$ such that, for all $f \in \cat_{x,y}$,
    \[ \id_y \circ f = f \circ \id_x = f.\] 
\end{enumerate}
\end{definition}
In Definition \ref{def:partialcat}, equality means that the left side is defined if and only if the right side is defined, and in this case, both sides are equal. In particular, the unitality axiom includes the requirement that $(f, \id_y) \in \cat_{x,y,y}$ and $(\id_x,f) \in \cat_{x,x,y}$.

A functor $F: \cat \to \cat'$ of partial categories is defined similarly to a functor of categories, with the additional specification that, if $f \circ g$ is undefined, then $F(f) \circ F(g)$ is undefined. We use $\catname{PCat}$ to denote the category of small partial categories.

The nerve construction for a category extends to partial categories in a straightforward way.

\begin{definition}
    The \emph{nerve} of a partial category $\cat$, denoted $N\cat_\bullet$, is given by 
    \[N\cat_n = \Hom_{\catname{PCat}}([n],\cat).\]
\end{definition}
The nerve of $\cat$ can be concretely described in low degrees by $N\cat_0 = \Ob(\cat)$, $N\cat_1 = \bigsqcup \cat_{x,y}$, and $N\cat_2 = \bigsqcup \cat_{x,y,z}$. More generally, for $n \geq 1$, we can identify $N\cat_n$ with the set of $n$-tuples $(f_1,f_2,\dots, f_n)$ that are fully composable (in the sense that $f_n \circ \cdots \circ f_1$ is defined), and the face and degeneracy maps are given by the same formulas as for the nerve of a category (see Section \ref{sec:nervecat}).

For the nerve of a partial category, the $n=2$ spine map is the inclusion map $\bigsqcup \cat_{x,y,z} \subseteq \bigsqcup \left(\cat_{x,y} \times \cat_{y,z}\right)$. It follows that $N\cat_\bullet$ is not Segal, except in the case where $\cat$ is actually a category.

\begin{proposition}\label{prop:partialcat}
    The nerve of a partial category is $2$-Segal.
\end{proposition}
\begin{proof}
    When $X_\bullet = N\cat_\bullet$ is the nerve of a partial category, the diagram \eqref{diag:subdivide2} is given by
\begin{equation*}
\begin{tikzcd}
{(f_1,\dots,f_n)} \arrow[d, maps to] \arrow[r, maps to] & {(f_1,\dots,f_i, f_j\circ \cdots \circ f_{i+1}, f_{j+1}, \dots, f_n)} \arrow[d, maps to] \\
{(f_{i+1},\dots,f_j)} \arrow[r, maps to]                & f_j\circ \cdots \circ f_{i+1}   
\end{tikzcd}
\end{equation*}
One can directly verify that this diagram is a pullback.
\end{proof}

Putting Propositions \ref{prop:segal} and \ref{prop:partialcat} together, we get
\begin{corollary}[\cite{Dyckerhoff-Kapranov:Higher}*{Proposition 2.3.4}]
    Every Segal set is $2$-Segal.
\end{corollary}

We emphasize that the converse to Proposition \ref{prop:partialcat} is not true, but there is a fairly simple condition (whose proof we leave as an exercise) that characterizes the $2$-Segal sets that are nerves of partial categories. 

\begin{proposition}\label{prop:partialcatnerve}
    A $2$-Segal set is the nerve of a partial category if and only if $T_2: X_2 \to X_1 \times_{X_0} X_1$ is injective.
\end{proposition}

\begin{example}[Partial monoids]\label{ex:partialmonoid}
A partial monoid \cite{Segal:Conf} is the special case of a partial category when $\Ob(\cat)$ has one element. The following are examples of partial monoids.
\begin{itemize}
    \item Monoids (and in particular, groups).
    \item For any positive integer $L$, the set $\mathcal{L} = \{1, x, x^2, \dots, x^L\}$ can be given a partial monoid structure where the operation is given by $x^i \circ x^j = x^{i+j}$ when $i+j \leq L$ and $x^i \circ x^j$ is undefined when $i+j > L$.
    \item For any set $S$, the power set $P(S)$ can be given a partial monoid structure where the operation is ``disjoint union'', i.e.\ $A \circ B = A \cup B$ when $A \cap B = \emptyset$ and $A \circ B$ is undefined when $A \cap B$ is nonempty.
\end{itemize}
\end{example}

\begin{example} [Quivers]\label{ex:quiver}
Consider a quiver where $V$ is the set of vertices and $E$ is the set of edges, with source and target maps $s,t: V \to E$. From this data, we can define a partial category $\mathcal{Q}$ as follows:
\begin{itemize}
    \item $\Ob(\mathcal{Q}) = V$.
    \item If $x \neq y$, then the set of morphisms from $x$ to $y$ is 
    \[ \mathcal{Q}_{x,y} = \{e \in E \suchthat s(e) = x, t(e) = y\}.\]
    If $x=y$, then
    \[ \mathcal{Q}_{x,x} = \{e \in E \suchthat s(e) = t(e) = x\} \cup \{\id_x\}.\]
    Thus the set of all morphisms in $\mathcal{Q}$ can be identified with $E \sqcup V$.
    \item The only compositions that are defined are the compositions with identity morphisms:
    \begin{align*}
        \id_y \circ e &= e \circ \id_x = e, & \id_x \circ \id_x = \id_x,
    \end{align*}
    where $e$ is an edge from $x$ to $y$.
\end{itemize}
The nerve of $\mathcal{Q}$ is the $1$-skeletal simplicial set $E \sqcup V \rightrightarrows V$. This example appears in \cite{Dyckerhoff-Kapranov:Higher}*{Example 3.1.1}.
\end{example}

\begin{example}[Cobordisms with restricted genus]
In \cite{BOORS}*{Example 2.2}, there is an example of a $2$-Segal set, denoted $\catname{2Cob}^{\leq g}_\bullet$, constructed using $2$-dimensional cobordisms with a restriction on the genus. This $2$-Segal set can be seen as the nerve of the partial category whose objects are $1$-dimensional closed manifolds (i.e.\ disjoint unions of circles), morphisms are diffeomorphism classes of cobordisms with genus $\leq g$, and where composition is given as usual by gluing of cobordisms, with the restriction that composition is undefined when gluing produces a cobordism with genus greater than $g$.
\end{example}

\begin{remark}
The inclusion of categories $\Cat \hookrightarrow \catname{PCat}$ has a left adjoint, giving a natural way to complete a partial category to a category with the same objects. Briefly, the completion of a partial category is obtained by taking formal compositions of morphisms, modulo the equivalence relation generated by the partially-defined composition. For details, see \cite{hermida-mateus1}*{Proposition 2.1}.

The nerve of a partial category (which is $2$-Segal) naturally sits inside of the nerve of its completion (which is $1$-Segal). For example, the completion of the partial monoid $\mathcal{L}$ in Example \ref{ex:partialmonoid} is $\{x^i \suchthat i \in \N \}\cong \N$, and the nerve of $\mathcal{L}$ is given by
\[ N\mathcal{L}_n = \left\{(x^{a_1},\dots,x^{a_n}) \in \mathcal{L}^n \suchthat \sum a_i \leq L\right\} \subset \{x^i \suchthat i \in \N \}^n \cong \N^n.\]

The completion of the partial category $\mathcal{Q}$ in Example \ref{ex:quiver} is the free category on the quiver, and the completion of $\catname{2Cob}^{\leq g}$ is the cobordism category $\catname{2Cob}$ (without any genus restriction).
\end{remark}

\subsection{\texorpdfstring{More examples of $2$-Segal sets}{More examples of 2-Segal sets}}
\label{sec:moreexamples}

There are many interesting examples of $2$-Segal sets that are not isomorphic to the nerve of a partial category; for instance, see the examples in \cite{bergner2024combinatorial} and \cite{GKT1} coming from graphs and trees. We also refer the reader to \cites{BOORS, CMS} and \cite{Dyckerhoff-Kapranov:Higher}*{Sections 3.1 and 3.2}.

Here, we give a couple of examples that can be described fairly explicitly. For both, we leave the verification of the simplicial identities and the $2$-Segal conditions as an exercise.

\begin{example}\label{ex:powerset}
    Let $A$ be a set, and let $P(A)$ denote the power set of $A$. We then form a simplicial set $X_\bullet$, where 
    \[X_n = \left\{(B_0, \dots, B_n) \in P(A)^{n+1} \suchthat \bigcup B_i = A\right\}.\]
    In particular, $X_0 = \{A\}$ only has one element. The face maps are given by
    \[ d_i^n(B_0, \dots, B_n) = (B_0, \dots, B_{i-1}, B_i \cup B_{i+1}, B_{i+2}, \dots, B_n)\]
    for $0 \leq i < n$, and 
    \[ d_n^n(B_0, \dots, B_n) = (B_0 \cup B_n, B_1, \dots, B_{n-1}),\]
    and the degeneracy maps are given by insertion of the empty set:
    \[ s_i^n(B_0, \dots, B_n) = (B_0, \dots, B_i, \emptyset, B_{i+1}, \dots, B_n).\]

To see that this simplicial set is not the nerve of a partial category (except when $A = \emptyset$), observe that
\begin{align*}
    T_2: X_2 &\to X_1 \times_{X_0} X_1, \\
    (B_0,B_1,B_2) &\mapsto ((B_0 \cup B_2, B_1), (B_0 \cup B_1, B_2)),
\end{align*}
is not injective, since, for example, $T_2(B,A,A) = ((A,A),(A,A))$ for all $B \in P(A)$. 
\end{example}

\begin{example} \label{ex:exponential}
For a natural number $m \geq 1$, we write $\underline{m} = \{1,\dots,m\}$, with $\underline{0} = \emptyset$. Let $X_n$ be the set whose elements are of the form $(m; S_1, \dots, S_n)$, where $S_1, \dots, S_n$ are (possibly empty) sets that form a partition of $\underline{m}$. In particular, $X_0 = \{0\}$ only has one element, and $X_1 = \{(m;\underline{m})\}$ is naturally isomorphic to $\N$. The face maps are given by nerve-like formulas:
\begin{align*}
    d_0^n(m; S_1, \dots, S_n) &= (m - |S_1|; S_2, \dots, S_n), \\
    d_i^n(m; S_1, \dots, S_n) &= (m; S_1, \dots, S_i \cup S_{i+1}, \dots, S_n), & 0 <i < n, \\
    d_n^n(m; S_1, \dots, S_n) &= (m - |S_n|; S_1, \dots, S_{n-1}),
\end{align*}
where in the formulas for $d_0^n$ and $d_n^n$ we have implicitly used the unique order-preserving isomorphisms $\underline{m-|S_1|} \cong \underline{m} \smallsetminus S_1$ and $\underline{m-|S_n|} \cong \underline{m} \smallsetminus S_n$. The degeneracy maps are given by insertion of the empty set. This example closely resembles the construction in \cite{BOORS}*{Example 2.3}, though it doesn't seem to precisely fit as a special case.
\end{example}

\subsection{\texorpdfstring{Graphical calculus for $2$-Segal sets}{Graphical calculus for 2-Segal sets}} \label{sec:gcalc}

Let $X_\bullet$ be a simplicial set. As we saw in Section \ref{sec:2segaldefinition}, there is a map \eqref{eqn:subdivide3} associated to each subdivision of the $(n+1)$-gon. These maps can be composed, since a polygon within a subdivision can itself be subdivided. Taken together, the subdivision maps and their compositions form a nice structure that we will now describe.

For each $n \geq 1$, consider the poset of subdivided $(n+1)$-gons, where the ordering is given by subordinate subdivisions. The first nontrivial cases, $n=3$ and $n=4$, are shown in Figures \ref{fig:assoc3} and \ref{fig:assoc4}, respectively.

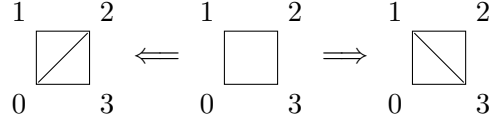
\begin{figure}[t]
\begin{center}
\begin{tikzpicture}[scale=0.5]
    \begin{scope} [shift={(5,0)}]
        \draw (-135:1) node[vertex] (0) {} node[anchor=north east] {0}
        -- (-45:1) node[vertex] (3) {} node[anchor=north west] {3}
        -- (45:1) node[vertex] (2) {} node[anchor=south west] {2}
        -- (135:1) node[vertex] (1) {} node[anchor=south east] {1}
        -- cycle;
    \end{scope}
    \node at (7.5,0) {$\Longrightarrow$};

        \begin{scope}[shift={(10,0)}]
        \draw (-135:1) node[vertex] (0) {} node[anchor=north east] {0}
        -- (-45:1) node[vertex] (3) {} node[anchor=north west] {3}
        -- (45:1) node[vertex] (2) {} node[anchor=south west] {2}
        -- (135:1) node[vertex] (1) {} node[anchor=south east] {1}
        -- cycle;
        \draw (3) -- (1);
    \end{scope}
    \node at (2.5,0) {$\Longleftarrow$};
    \begin{scope}
        \draw (-135:1) node[vertex] (0) {} node[anchor=north east] {0}
        -- (-45:1) node[vertex] (3) {} node[anchor=north west] {3}
        -- (45:1) node[vertex] (2) {} node[anchor=south west] {2}
        -- (135:1) node[vertex] (1) {} node[anchor=south east] {1}
        -- cycle;
        \draw (2) -- (0);
    \end{scope}
\end{tikzpicture}
\end{center}
\caption{The poset of subdivided squares.}
\label{fig:assoc3}
\end{figure}
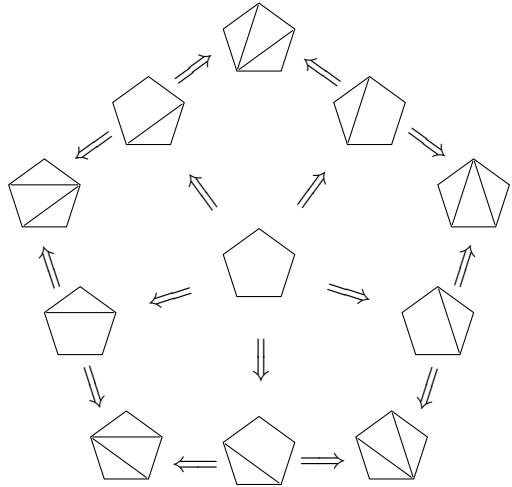
\begin{figure}[t]
\begin{center}
    \begin{tikzpicture}[scale=0.5]
            \begin{scope}
        \draw (-126:1) node[vertex] (0) {} 
        -- (-54: 1) node[vertex] (4) {} 
        -- (18:1) node[vertex] (3) {} 
        -- (90:1) node[vertex] (2) {} 
        -- (162:1) node[vertex] (1) {} 
        -- cycle;
    \end{scope}
    \begin{scope}[shift={(54:5)}]
        \draw (-126:1) node[vertex] (0) {} 
        -- (-54: 1) node[vertex] (4) {} 
        -- (18:1) node[vertex] (3) {} 
        -- (90:1) node[vertex] (2) {} 
        -- (162:1) node[vertex] (1) {} 
        -- cycle;
        \draw (0) -- (2);
    \end{scope}
        \begin{scope}[shift={(90:6)}]
        \draw (-126:1) node[vertex] (0) {} 
        -- (-54: 1) node[vertex] (4) {} 
        -- (18:1) node[vertex] (3) {} 
        -- (90:1) node[vertex] (2) {} 
        -- (162:1) node[vertex] (1) {} 
        -- cycle;
        \draw (0) -- (2);
        \draw (0) -- (3);
    \end{scope}
        \begin{scope}[shift={(126:5)}]
        \draw (-126:1) node[vertex] (0) {} 
        -- (-54: 1) node[vertex] (4) {} 
        -- (18:1) node[vertex] (3) {} 
        -- (90:1) node[vertex] (2) {} 
        -- (162:1) node[vertex] (1) {} 
        -- cycle;
        \draw (0) -- (3);
    \end{scope}
        \begin{scope}[shift={(162:6)}]
        \draw (-126:1) node[vertex] (0) {} 
        -- (-54: 1) node[vertex] (4) {} 
        -- (18:1) node[vertex] (3) {} 
        -- (90:1) node[vertex] (2) {} 
        -- (162:1) node[vertex] (1) {} 
        -- cycle;
        \draw (0) -- (3);
        \draw (1) -- (3);
    \end{scope}
        \begin{scope}[shift={(198:5)}]
        \draw (-126:1) node[vertex] (0) {} 
        -- (-54: 1) node[vertex] (4) {} 
        -- (18:1) node[vertex] (3) {} 
        -- (90:1) node[vertex] (2) {} 
        -- (162:1) node[vertex] (1) {} 
        -- cycle;
        \draw (1) -- (3);
    \end{scope}
         \begin{scope}[shift={(234:6)}]
        \draw (-126:1) node[vertex] (0) {} 
        -- (-54: 1) node[vertex] (4) {} 
        -- (18:1) node[vertex] (3) {} 
        -- (90:1) node[vertex] (2) {} 
        -- (162:1) node[vertex] (1) {} 
        -- cycle;
        \draw (1) -- (3);
        \draw (1) -- (4);
    \end{scope}
        \begin{scope}[shift={(270:5)}]
        \draw (-126:1) node[vertex] (0) {} 
        -- (-54: 1) node[vertex] (4) {} 
        -- (18:1) node[vertex] (3) {} 
        -- (90:1) node[vertex] (2) {} 
        -- (162:1) node[vertex] (1) {} 
        -- cycle;
        \draw (1) -- (4);
    \end{scope}
        \begin{scope}[shift={(306:6)}]
        \draw (-126:1) node[vertex] (0) {} 
        -- (-54: 1) node[vertex] (4) {} 
        -- (18:1) node[vertex] (3) {} 
        -- (90:1) node[vertex] (2) {} 
        -- (162:1) node[vertex] (1) {} 
        -- cycle;
        \draw (1) -- (4);
        \draw (2) -- (4);
    \end{scope}
        \begin{scope}[shift={(342:5)}]
        \draw (-126:1) node[vertex] (0) {} 
        -- (-54: 1) node[vertex] (4) {} 
        -- (18:1) node[vertex] (3) {} 
        -- (90:1) node[vertex] (2) {} 
        -- (162:1) node[vertex] (1) {} 
        -- cycle;
        \draw (2) -- (4);
    \end{scope}
        \begin{scope}[shift={(18:6)}]
        \draw (-126:1) node[vertex] (0) {} 
        -- (-54: 1) node[vertex] (4) {} 
        -- (18:1) node[vertex] (3) {} 
        -- (90:1) node[vertex] (2) {} 
        -- (162:1) node[vertex] (1) {} 
        -- cycle;
        \draw (2) -- (4);
        \draw (0) -- (2);
    \end{scope}
    \node[rotate=342] at (342:2.5) {$\Longrightarrow$};
    \node[rotate=54] at (54:2.5) {$\Longrightarrow$};
    \node[rotate=126] at (126:2.5) {$\Longrightarrow$};
    \node[rotate=198] at (198:2.5) {$\Longrightarrow$};
    \node[rotate=270] at (270:2.5) {$\Longrightarrow$};
    \node[rotate=72] at (0:5.5) {$\Longrightarrow$};    
    \node[rotate=144] at (72:5.5) {$\Longrightarrow$};
    \node[rotate=216] at (144:5.5) {$\Longrightarrow$};
    \node[rotate=288] at (216:5.5) {$\Longrightarrow$};
    \node[rotate=0] at (288:5.5) {$\Longrightarrow$};
    \node[rotate=324] at (36:5.5) {$\Longrightarrow$};    
    \node[rotate=36] at (108:5.5) {$\Longrightarrow$};
    \node[rotate=108] at (180:5.5) {$\Longrightarrow$};
    \node[rotate=180] at (252:5.5) {$\Longrightarrow$};
    \node[rotate=252] at (324:5.5) {$\Longrightarrow$};
    \end{tikzpicture}
\end{center}
\caption{The poset of subdivided pentagons.}
\label{fig:assoc4} 
\end{figure}

Collectively, the maps \eqref{eqn:subdivide3} can be viewed as giving, for each $n \geq 2$, a functor from the poset category of subdivided $(n+1)$-gons to the category of sets. In \cite{CMS}, this functor was called the \emph{$2$-Segal functor}. For $n=3$, the $2$-Segal functor takes the poset in Figure \ref{fig:assoc3} to the diagram
\begin{equation}\label{eqn:tacos}
X_2 \times_{X_1} X_2 \xlongleftarrow{\mathcal{T}_{13}} X_3 \xlongrightarrow{\mathcal{T}_{02}} X_2 \times_{X_1} X_2,
\end{equation}
where
\begin{align*}
\mathcal{T}_{13}&: \psi \mapsto (d_1 \psi, d_3 \psi),\\
\mathcal{T}_{02}&: \psi \mapsto (d_0 \psi, d_2 \psi).
\end{align*}
The first step toward a graphical calculus is the observation that we can represent (for example) the maps in \eqref{eqn:tacos} by the pictures in Figure \ref{fig:assoc3}. 

We stress that there are two different fiber products\footnote{These are the $(13)$- and $(02)$-\emph{taco spaces}, as named in \cite{ContrerasMehtaSpan}.} in \eqref{eqn:tacos}; the one on the left is
\[ \{(\omega_1, \omega_3) \in (X_2)^2 \suchthat d_2 \omega_1 = d_1 \omega_3\},\]
and the one on the right is
\[ \{(\omega_0, \omega_2) \in (X_2)^2 \suchthat d_1 \omega_0 = d_0 \omega_2\}.\]
This distinction is absent in the notation of \eqref{eqn:tacos}, whereas it can be inferred from the configuration of the subdivisions in Figure \ref{fig:assoc3}.

As an immediate consequence of the definitions, we have
\begin{proposition}\label{prop:2segalfunctor}
    A simplicial set $X_\bullet$ is $2$-Segal if and only if, for all $n \geq 2$, the image of the associated $2$-Segal functor consists of isomorphisms.
\end{proposition}

For the remainder of this section, we will assume that $X_\bullet$ is $2$-Segal. Then, Proposition \ref{prop:2segalfunctor} justifies the use of graphical representations as follows:
\begin{itemize}
    \item An unsubdivided $(n+1)$-gon represents $X_n$.
    \item A subdivided $(n+1)$-gon represents the target space for the associated map \eqref{eqn:subdivide3}.
    \item The $2$-Segal functor provides canonical isomorphisms between the spaces represented by any two subdivisions of the $(n+1)$-gon.
\end{itemize}
Frequently, we will label a polygon to indicate a specific element of the set that the polygon represents.

The true value of these graphical representations comes from the way that they interact with face and degeneracy maps. We first consider the face maps. By direct comparison of the formulas, one can see that the coface maps $\delta_i^n: [n-1] \to [n]$ are actually special cases of the maps introduced in Definition \ref{def:PQ}. Specifically,
\begin{equation}\label{eqn:facepq}
  \delta_i^n = \begin{cases}
    P_{1n}^n, & i=0,\\
    Q_{(i-1)(i+1)}^n, & 0 < i < n,\\
    P_{0(n-1)}^n, & i=n.
\end{cases}  
\end{equation}
Thus the face maps can be recognized as components of the corresponding maps in \eqref{diag:subdivide2}. 

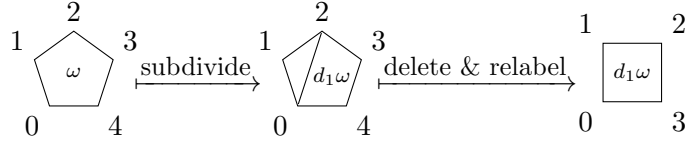
\begin{figure}[t]
\begin{center}
\begin{tikzpicture}[scale=0.55]
    \begin{scope}
        \draw (-126:1) node[vertex] (0) {} node[anchor=north east] {0}
        -- (-54: 1) node[vertex] (4) {} node[anchor=north west] {4}
        -- (18:1) node[vertex] (3) {} node[anchor=south west] {3}
        -- (90:1) node[vertex] (2) {} node[anchor=south] {2}
        -- (162:1) node[vertex] (1) {} node[anchor=south east] {1}
        -- cycle;
        \node at (0,0) {${\scriptstyle \omega}$};
    \end{scope}
    \node at (3,0) {$\xmapsto{\mbox{\small{subdivide}}}$};
    \begin{scope}[shift={(6,0)}]
        \draw (-126:1) node[vertex] (0) {} node[anchor=north east] {0}
        -- (-54: 1) node[vertex] (4) {} node[anchor=north west] {4}
        -- (18:1) node[vertex] (3) {} node[anchor=south west] {3}
        -- (90:1) node[vertex] (2) {} node[anchor=south] {2}
        -- (162:1) node[vertex] (1) {} node[anchor=south east] {1}
        -- cycle;
        \node at (-18:.25) {${\scriptstyle d_1 \omega}$};
        \draw (0) -- (2);
    \end{scope}
    \node at (9.75,0) {$\xmapsto{\mbox{\small{delete \& relabel}}}$};
    \begin{scope}[shift={(13.5,0)}]
        \draw (-135:1) node[vertex] (0) {} node[anchor=north east] {0}
        -- (-45:1) node[vertex] (3) {} node[anchor=north west] {3}
        -- (45:1) node[vertex] (2) {} node[anchor=south west] {2}
        -- (135:1) node[vertex] (1) {} node[anchor=south east] {1}
        -- cycle;
        \node at (0,0) {${\scriptstyle d_1\omega}$};
    \end{scope}
\end{tikzpicture}
\end{center}
    \caption{Graphical calculus for the face map $d_1^4: X_4 \to X_3$.}
    \label{fig:facemaps}
\end{figure}

In terms of graphical representations, all three cases in \eqref{eqn:facepq} follow a common rule. In the $(n+1)$-gon, draw the edge connecting the two vertices that are adjacent to vertex $i$. This edge subdivides the $(n+1)$-gon into an $n$-gon and a triangle, and the face map can be seen as restricting to the $n$-gon component. In other words, we can visualize $d_i^n$ as the process of deleting the triangle formed by vertex $i$ and its two adjacent vertices; see Figure \ref{fig:facemaps}.

\begin{figure}[t]
    \begin{center}
    \begin{tikzpicture}[scale=0.5]
            \begin{scope}
        \draw (240:1) node[vertex] (0) {} node[anchor=north east] {0}
        -- (300:1) node[vertex] (5) {} node[anchor=north west] {5}
        -- (0:1) node[vertex] (4) {} node[anchor=west] {4}
        -- (60:1) node[vertex] (3) {} node[anchor=south west] {3}
        -- (120:1) node[vertex] (2) {} node[anchor=south east] {2}
        -- (180:1) node[vertex] (1) {} node[anchor=east] {1}        
        -- cycle;
        \draw (0) -- (2);
        \draw (2) -- (4);
    \end{scope}
    \begin{scope}[shift={(8,0)}]
        \draw (240:1) node[vertex] (0) {} node[anchor=north east] {0}
        -- (300:1) node[vertex] (5) {} node[anchor=north west] {5}
        -- (0:1) node[vertex] (4) {} node[anchor=west] {4}
        -- (60:1) node[vertex] (3) {} node[anchor=south west] {3}
        -- (120:1) node[vertex] (2) {} node[anchor=south east] {2}
        -- (180:1) node[vertex] (1) {} node[anchor=east] {1}        
        -- cycle;
        \draw (1) -- (4);
        \draw (2) -- (4);
    \end{scope}
\node at (10.5,0) {$\overset{\sim}\Longrightarrow$};
        \begin{scope}[shift={(13,0)}]
        \draw (240:1) node[vertex] (0) {} node[anchor=north east] {0}
        -- (300:1) node[vertex] (5) {} node[anchor=north west] {5}
        -- (0:1) node[vertex] (4) {} node[anchor=west] {4}
        -- (60:1) node[vertex] (3) {} node[anchor=south west] {3}
        -- (120:1) node[vertex] (2) {} node[anchor=south east] {2}
        -- (180:1) node[vertex] (1) {} node[anchor=east] {1}        
        -- cycle;
        \draw (1) -- (4);
        \draw (1) -- (3);
    \end{scope}
    \end{tikzpicture}
    \end{center}
    \caption{On the left, visualization of the identity $d_1^4 d_3^5 = d_2^4 d_1^5$. Both sides have the effect of deleting the same two triangles. On the right, visualization of the identity $d_2^4 d_3^5 = d_2^4 d_2^5$. Both sides have the effect of deleting the same quadrilateral.
    }
    \label{fig:faceface}
\end{figure}
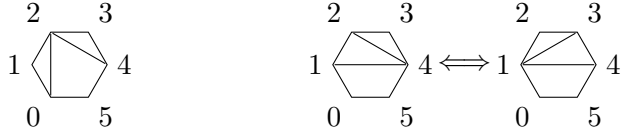

The graphical calculus provides a nice way to visualize the face map identities \eqref{eqn:face1}. For example, if $i < j-1$, then the triangles that are deleted by $d_i$ and $d_j$ are non-overlapping, so it follows that the operations commute (up to an index shift), giving the identity $d_i d_j = d_{j-1}d_i$. In the case where $i = j-1$, one can see that $d_i d_j = d_{j-1}d_i$ because both sides have the effect of deleting the same quadrilateral; see Figure \ref{fig:faceface}.

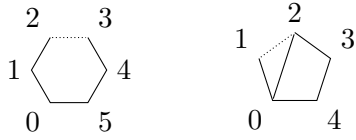
\begin{figure}[t]
    \begin{center}
    \begin{tikzpicture}[scale=0.5]
            \begin{scope}
        \draw (120:1) node[vertex] (2) {} node[anchor=south east] {2}
        -- (180:1) node[vertex] (1) {} node[anchor=east] {1}
        -- (240:1) node[vertex] (0) {} node[anchor=north east] {0}
        -- (300:1) node[vertex] (5) {} node[anchor=north west] {5}
        -- (0:1) node[vertex] (4) {} node[anchor=west] {4}
        -- (60:1) node[vertex] (3) {} node[anchor=south west] {3};
        \draw[densely dotted] (2) -- (3);
    \end{scope}
    \begin{scope}[shift={(6,0)}]
        \draw 
        (162:1) node[vertex] (1) {} node[anchor=south east] {1}
        --(-126:1) node[vertex] (0) {} node[anchor=north east] {0}
        -- (-54: 1) node[vertex] (4) {} node[anchor=north west] {4}
        -- (18:1) node[vertex] (3) {} node[anchor=south west] {3}
        -- (90:1) node[vertex] (2) {} node[anchor=south] {2};
        \draw[densely dotted] (1) -- (2);
        \draw (0) -- (2);
    \end{scope}
    \end{tikzpicture}
    \end{center}
    \caption{The figure on the left represents $\{\omega \in X_5 \suchthat e_3 \omega \in s_0(X_0)\}$. The figure on the right represents $\{(\xi,\psi) \in X_2 \times X_3 \suchthat e_\out \xi = e_1 \psi,\; e_2 \xi \in s_0(X_0)\}$.}
    \label{fig:degenerate}
\end{figure}

We now turn to the degeneracy maps. For these, it will be useful to use a dotted edge to represent the set $s_0(X_0) \subseteq X_1$ of degenerate $1$-simplices. An $(n+1)$-gon with a dotted edge represents the subset of $X_n$ for which the corresponding edge is degenerate. Although in principle any edge can be dotted, we will only make use of diagrams where external edges corresponding to vertebrae are dotted. A couple of examples are illustrated in Figure \ref{fig:degenerate}.

\begin{figure}[t]
    \begin{center}
    \begin{tikzpicture}[scale=0.5]
    \begin{scope}
        \draw (120:1) node[vertex] (2) {} node[anchor=south east] {2}
        -- (180:1) node[vertex] (1) {} node[anchor=east] {1}
        -- (240:1) node[vertex] (0) {} node[anchor=north east] {0}
        -- (300:1) node[vertex] (5) {} node[anchor=north west] {5}
        -- (0:1) node[vertex] (4) {} node[anchor=west] {4}
        -- (60:1) node[vertex] (3) {} node[anchor=south west] {3};
        \draw[densely dotted] (2) -- (3);
        \node at (0,0) {${\scriptstyle \omega}$};
    \end{scope}
    \node at (2.5,0) {$\overset{\sim}\Longrightarrow$};
    \begin{scope}[shift={(5,0)}]
        \draw (120:1) node[vertex] (2) {} node[anchor=south east] {2}
        -- (180:1) node[vertex] (1) {} node[anchor=east] {1}
        -- (240:1) node[vertex] (0) {} node[anchor=north east] {0}
        -- (300:1) node[vertex] (5) {} node[anchor=north west] {5}
        -- (0:1) node[vertex] (4) {} node[anchor=west] {4}
        -- (60:1) node[vertex] (3) {} node[anchor=south west] {3};
        \draw[densely dotted] (2) -- (3);
        \draw (1) -- (3);
        \node at (300:.2) {${\scriptstyle \eta}$};
    \end{scope}    
    \node at (7.5,0) {$\overset{\sim}\Longrightarrow$};
    \begin{scope}[shift={(10,0)}]
        \draw (120:1) node[vertex] (2) {} node[anchor=south east] {2}
        -- (180:1) node[vertex] (1) {} node[anchor=east] {1}
        -- (240:1) node[vertex] (0) {} node[anchor=north east] {0}
        -- (300:1) node[vertex] (5) {} node[anchor=north west] {5}
        -- (0:1) node[vertex] (4) {} node[anchor=west] {4}
        -- (60:1) node[vertex] (3) {} node[anchor=south west] {3};
        \draw[densely dotted] (2) -- (3);
        \draw (2) -- (4);
        \node at (240:.2) {${\scriptstyle \eta}$};
    \end{scope} 
    \end{tikzpicture}
    \end{center}
    \caption{Since $e_3 \omega$ is degenerate, we have that $\omega = s_2 \eta$, where $\eta = d_2 \omega = d_3 \omega$. This also illustrates that there are two ways to produce $s_2 \eta$ by appending a degenerate $2$-simplex so that the result has degenerate $3$rd vertebra.}
    \label{fig:higherdegenerate}
\end{figure}
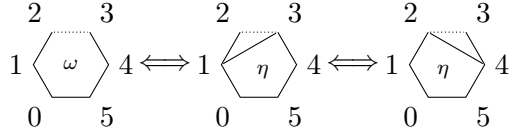

Recall (see Remark \ref{remark:unitality}) that the diagram
\begin{equation*}
\begin{tikzcd}
X_{n-1} \arrow[r, "s_i"] \arrow[d, "v_i"'] & X_n \arrow[d, "e_{i+1}"] \\
X_0 \arrow[r, "s_0"] & X_1
\end{tikzcd}
\end{equation*}
is a pullback by the unitality property of $2$-Segal sets. Since degeneracy maps are injective, it follows that the set represented by the $(n+1)$-gon with dotted $(i+1)$th vertebra is
\begin{equation*}
    \{\omega \in X_n \suchthat e_{i+1} \omega \in s_0(X_0)\} = s_i(X_{n-1}).
\end{equation*}
Furthermore, if $\omega = s_i \eta$, then it follows from the simplicial identities that $\eta = d_i \omega = d_{i+1}\omega$. This relationship allows us to visualize the degeneracy map $s_i^n$ as the process of appending a degenerate $2$-simplex so that the resulting $(n+2)$-gon has a degenerate $(i+1)$th vertebra (see Figure \ref{fig:higherdegenerate}).

We end this section with a summary of the graphical calculus for $2$-Segal sets:
\begin{itemize}
    \item The face map $d_i$ is obtained by deleting the $2$-simplex formed by the vertex $i$ and its two adjacent vertices.
    \item The degeneracy map $s_i$ is obtained by appending a degenerate $2$-simplex so that the result has a degenerate $(i+1)$th vertebra.
\end{itemize}

\section{Pseudomonoids in the bicategory of spans}
\label{sec:span}

In this section, we briefly review the bicategory $\Span$, its monoidal structure, and the notion of pseudomonoid. For more details on $\Span$, we refer the reader to \cite{CKWW}, which gives a general construction of the monoidal bicategory of spans in any category with finite limits, and to \cite{stay}, where a clear and complete definition of monoidal bicategories is given, with bicategories of spans as the main example. For pseudomonoids, the original reference is \cite{DayStreet}, and we also mention \cite{Verdon}, which uses string diagrams similar to those presented here.

We note that, while $\Span$ has a natural \emph{symmetric} monoidal structure (and this is discussed in the above references), the symmetric structure does not play any role in the present paper. In \cite{CMS}, the symmetric structure is used to define \emph{commutative} pseudomonoids in $\Span$, which are then characterized in terms of $2$-Segal sets with additional structure.

\subsection{\texorpdfstring{The bicategory $\Span$}{The bicategory Span}}\label{sec:bicatspan}

\begin{definition}
    
The bicategory $\Span$ is given by the following data:
\begin{itemize}
    \item objects are sets;
    \item a morphism from a set $X$ to a set $Y$ is a span $X \xleftarrow{f_1} A \xrightarrow{f_2} Y$, which we will sometimes also denote as $f: X \leftarrow A \rightarrow Y$;
    \item a $2$-morphism from $f: X \xleftarrow{f_1} A \xrightarrow{f_2} Y$ to $g: X \xleftarrow{g_1} B \xrightarrow{g_2} Y$ is a map of spans, i.e.\ a map $h: A \to B$ such that 
\[\begin{tikzcd}[sep=small]
  & A \arrow[ld, "f_1"'] \arrow[rd, "f_2"] \arrow[dd, "h"] &   \\
X &                                                        & Y \\
  & B \arrow[lu, "g_1"] \arrow[ru, "g_2"']                 &  
\end{tikzcd}\]
commutes; such a $2$-morphism will usually be denoted by a double-arrow, as in $f \xRightarrow{h} g$;
    \item the identity morphism is the span $\id_X: X \xleftarrow{\id} X \xrightarrow{\id} X$;
    \item the composition of $f: X \xleftarrow{f_1} A \xrightarrow{f_2} Y$ with $g: Y \xleftarrow{g_1} B \xrightarrow{g_2} Z$ is formed by taking the pullback
    \[\begin{tikzcd}[sep=small]
  &                                        & A \times_Y B \arrow[ld] \arrow[rd] &                                        &   \\
  & A \arrow[ld, "f_1"'] \arrow[rd, "f_2"] &                                    & B \arrow[ld, "g_1"'] \arrow[rd, "g_2"] &   \\
X &                                        & Y                                  &                                        & Z
\end{tikzcd}\]
giving the span $g \circ f: X \leftarrow A \times_Y B \rightarrow Z$.
\end{itemize}

\end{definition}

The remaining bicategorical data (horizontal composition of $2$-morphisms, associator, and unitors) arise naturally from the universal property of pullbacks. For example, the associator is given by the canonical isomorphism 
\begin{equation}\label{eqn:spanassociator}
(A \times_X B) \times_Y C \cong A \times_X (B \times_Y C),
\end{equation}
and the unitors are given by the canonical isomorphisms 
\[X \times_X A \cong A \cong A \times_X X.\]

Given sets $X,Y$, it can sometimes be convenient to consider the category $\Hom_\Span(X,Y)$, whose objects are spans $X \leftarrow A \rightarrow Y$, and whose morphisms are maps of spans. (Equivalently, $\Hom_\Span(X,Y)$ can be viewed as the slice category $\set/(X \times Y)$.) Thus, the Hom-categories $\Hom_\Span(X,Y)$ combine the morphisms and $2$-morphisms into a single entity. We can then view composition as a functor $\Hom_\Span(Y,Z) \times \Hom_\Span(X,Y) \to \Hom_\Span(X,Z)$.

\subsection{\texorpdfstring{The monoidal structure of $\Span$}{The monoidal structure of Span}}

The monoidal structure on $\Span$ is given by the Cartesian product, with the monoidal unit (some fixed choice of) a $1$-point set, which we will denote $\pt$. Again, all of the higher structures arise naturally from universal properties, but we highlight one particular structure map that will appear in the definition of pseudomonoid.

Part of the data of a monoidal bicategory is what Stay \cite{stay} calls the \emph{tensorator}, which is a $2$-isomorphism controlling the failure of the monoidal product to commute with composition. We will now describe the tensorator in $\Span$. Suppose we have two pairs of composable spans 
\begin{align*}
X \leftarrow A \rightarrow &X' \leftarrow A' \rightarrow X'', \\
Y \leftarrow B \rightarrow &Y' \leftarrow B' \rightarrow Y''.     
\end{align*}
We could take the Cartesian product first, yielding
\[ X \times Y \leftarrow A \times B \rightarrow X' \times Y' \leftarrow A' \times B' \rightarrow X'' \times Y'',\]
and then compose, yielding a span
\begin{equation}\label{eqn:tensorator1}
X \times Y \leftarrow (A \times B) \times_{(X' \times Y')} (A' \times B') \rightarrow X'' \times Y''.    
\end{equation} 
Alternatively, we could do the compositions first, yielding spans 
\begin{align*}
X \leftarrow &A \times_{X'} A' \rightarrow X'', \\
Y \leftarrow &B \times_{Y'} B' \rightarrow Y'',
\end{align*}
and then take the Cartesian product, yielding a span
\begin{equation}\label{eqn:tensorator2}
    X \times Y \leftarrow (A \times_{X'} A') \times (B \times_{Y'} B') \rightarrow X'' \times Y''.
\end{equation} 
The tensorator is the canonical isomorphism between the spans \eqref{eqn:tensorator1} and \eqref{eqn:tensorator2}.

As a special case, if we have spans $f: X \leftarrow A \rightarrow X'$ and $g: Y \leftarrow B \rightarrow Y'$, the tensorator and unitors give 2-isomorphisms 
\[(\id_{X'} \times g) \circ (f \times \id_{Y}) \overset{\sim}\Longrightarrow (\id_{X'} \circ f) \times (g \circ \id_Y) \overset{\sim}\Longrightarrow f \times g\]
and
\[(f \times \id_{Y'}) \circ (\id_X \times g) \overset{\sim}\Longrightarrow (f \circ \id_X) \times (\id_{Y'} \circ g) \overset{\sim}\Longrightarrow f \times g\]
in $\Hom_\Span(X \times Y, X' \times Y')$, given by the canonical isomorphisms
    \[ (A \times Y) \times_{(X' \times Y)} (X' \times B) \cong  (A \times_{X'} X') \times (Y \times_Y B) \cong A \times B\]
    and
    \[ (X \times B) \times_{(X \times Y')} (A \times Y') \cong (X \times_X A) \times (B \times_{Y'} Y') \cong A \times B.\]
    Composing these isomorphisms gives a $2$-isomorphism
    \begin{equation}\label{eqn:slide}
    (\id_{X'} \times g) \circ (f \times \id_Y) \xRightarrow{c_{f,g}} (f \times \id_{Y'}) \circ (\id_X \times g).
    \end{equation}
The map $c_{f,g}$ is called the \emph{slide move}, a name which will become more visually intuitive in the following section.

\subsection{String diagrams}

When working within a monoidal category, it can be convenient to use string diagrams, as described in \cite{joyal-street}. In a monoidal bicategory, string diagrams can still be used, with certain caveats (see Remark \ref{remark:caveat} below)\footnote{For bicategories, it is common to use diagrams that are dual to the ones we describe here. In this dual picture, objects, morphisms, and $2$-morphisms are represented by planar regions, lines, and points, respectively. However, monoidal structures are difficult to visualize from this perspective, requiring $3$-dimensional diagrams.}. We provide a brief review here. Because we are solely interested in $\Span$, we will use the notation $\times$ for the monoidal product and $\pt$ for the monoidal unit, even though the diagrams are equally valid in an arbitary monoidal bicategory.

A morphism $f$ from $X$ to $Y$ is drawn as 
		 \stringdiagram{
		 	\path (0,2.5) node {$X$};
                \halfidentity{0}{2}
		 	\morphism{0}{0}{f}
                \halfidentity{0}{-1}
		 	\path (0,-2.5) node {$Y$};
		}
In particular, we use the convention that diagrams are read from top-to-bottom\footnote{We warn the reader that this convention is not standard. It is common for people to draw diagrams bottom-to-top or left-to-right.}. When they are clear from context (as will be the case in most of this paper), the labels $X$ and $Y$ on the strings will be omitted.

Given composable morphisms $f$ and $g$, the diagram
		\stringdiagram{
		\morphism{0}{0}{f}
		\morphism{0}{-2}{g}
		}
represents the composition $g\circ f$. Monoidal products are shown by drawing diagrams in parallel, and the absence of strings represents the monoidal unit $\pt$. For example, the diagrams
		\stringdiagram
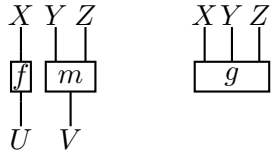
{
        \begin{scope}
		\halfidentity{0}{2}
		\morphism{0}{0}{f}
		\halfidentity{0}{-1} 
		 \draw (1,0.6) rectangle (3,-0.6);
		 \path (2,0) node {$m$}; 
		 \draw (1.4,0.6) -- (1.4,2);
		 \draw (2.6,0.6) -- (2.6,2);
		 \draw (2,-0.6) -- (2,-2); 
		 \path (0,2.5) node {$X$}; 
		 \path (1.4,2.5) node {$Y$}; 
		 \path (2.6,2.5) node {$Z$}; 
		 \path (0,-2.5) node {$U$};
		 \path (2,-2.5) node {$V$};  
         \end{scope}
            \begin{scope}[shift={(7,0)}]
          \draw (0,0.6) rectangle (3,-0.6);
		 \path (1.5,0) node {$g$}; 
          \draw (0.4,0.6) -- (0.4,2);
		 \draw (1.5,0.6) -- (1.5,2);
		 \draw (2.6,0.6) -- (2.6,2); 
          \path (0.4,2.5) node {$X$};
		 \path (1.5,2.5) node {$Y$}; 
		 \path (2.6,2.5) node {$Z$}; 
         \end{scope}
			}
		represent morphisms $f\times m :X\times (Y\times Z)\to U\times V$ and $g: X \times Y \times Z \to \pt$.

\begin{figure}[t]
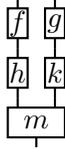

    \begin{center}
        		\stringdiagram{
		\morphism{0}{0}{f}
            \morphism{1.5}{0}{g}
            \morphism{0}{-2}{h}
            \morphism{1.5}{-2}{k}
		 \draw (-.4,-3.4) rectangle (1.9,-4.6);
		 \path (0.75,-4) node {$m$}; 
		 \draw (0,-3.4) -- (0,-3);
		 \draw (1.5,-3.4) -- (1.5,-3);
		 \draw (.75,-4.6) -- (.75,-5);  
			}
    \end{center}
    \caption{This diagram could represent $(m \circ (h \times k)) \circ (f \times g)$ or $m \circ ((h \times k) \circ(f \times g))$ or $m \circ ((f \circ h) \times (g \circ k))$. Insertion of identity morphisms or monoidal units is also possible, e.g.\ $(m \circ (\id \times \id) \circ (h \times k)) \circ (f \times g)$.}
    \label{fig:order}
\end{figure}

\begin{remark}\label{remark:caveat}
    An important caveat when using string diagrams in a monoidal \emph{bicategory} is that there can be ambiguities coming from the fact that associativity and unitality (for both composition and the monoidal product), as well as compatibility of composition with the monoidal product, only hold up to $2$-isomorphism. As a result, a given diagram can have many possible interpretations, depending on the order of operations and the possibility of inserting identity morphisms or monoidal units; see Figure \ref{fig:order}. However, this difficulty can be resolved by the coherence theorem for monoidal bicategories (see \cite{gurski:coherence} and the references within), which implies that any two interpretations of a diagram are related by a \emph{canonical} $2$-isomorphism\footnote{In particular, in $\Span$, the associators, unitors, etc.\ (e.g.\ \eqref{eqn:spanassociator}) are \emph{so} canonical that it's common to write them as equalities. For example, it's common for people to write $(X \times Y) \times Z = X \times (Y \times Z)$ as a shorthand for the fact that the two sides are canonically isomorphic.}. In practice, this means that, for any string diagram, one should assume that there is some chosen order of operations, but it is unnecessary to explicitly specify the choice.
\end{remark}
    
Finally, we note that string diagrams give us a nice visualization of the slide move \eqref{eqn:slide}, which explains the terminology:
    \stringdiagram{
    \morphism{0}{1}{f}
    \identity{1}{1}
    \identity{0}{-1}
    \morphism{1}{-1}{g}
    \nattrans{3}{0}{c_{f,g}}
    \morphism{5}{-1}{f}
    \identity{6}{-1}
    \identity{5}{1}
    \morphism{6}{1}{g}
    }

\subsection{Projection maps}\label{sec:projection}

There is an operation that is specific to $\Span$ that will be useful for us later. It comes from the observation that, in $\Span$, composition is given by pullback (in the category of sets), and the monoidal operation is the Cartesian product, which is also a pullback in the category of sets. Thus, in both cases, there are natural maps from the set forming the apex of the composite span to the set forming the apex of each component.

To be more explicit, consider spans
\begin{align}
    f&: X \leftarrow A \rightarrow Y, & g&: Y \leftarrow B \rightarrow Z. \label{eqn:spanscompose}
\end{align}
Then we can form the composition
\begin{align*}
    g \circ f &: X \leftarrow A \times_Y B \rightarrow Z.
\end{align*}
In this situation, we refer to the natural maps $A \times_Y B \to A$ and $A \times_Y B \to B$ as the \emph{projection maps} from $g \circ f$ to the components $f$ and $g$, respectively. 

Similarly, given spans
\begin{align}
    f&: X \leftarrow A \rightarrow Y, & h&: U \leftarrow C \rightarrow V, \label{eqn:spanscompose2}
\end{align}
we can form the product
\begin{align*}
    f \times h&: X \times U \leftarrow A \times C \rightarrow Y \times V,
\end{align*}
and we refer to the natural maps $A \times C \to A$ and $A \times C \to C$ as the projection maps from $f \times h$ to the components $f$ and $h$, respectively.

\begin{figure}[t]
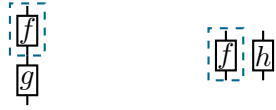

\begin{center}
  		\stringdiagram{
        \begin{scope}
		\morphism{0}{0}{f}
		\morphism{0}{-2}{g}
            \draw[dashed, Gray] (-.7,-1) rectangle (.7,1.1);
        \end{scope}
        \begin{scope}[shift={(8,-1)}]
		\morphism{0}{0}{f}
		\morphism{1.5}{0}{h}
            \draw[dashed, Gray] (-.7,-1) rectangle (.7,1.1);
        \end{scope}
		}    
\end{center}
 \caption{If $f,g,h$ are as in \eqref{eqn:spanscompose} and \eqref{eqn:spanscompose2}, the diagram on the left represents the projection map $A \times_Y B \to A$, and the diagram on the right represents the projection map $A\times C \to A$.}
    \label{fig:subdiagram}
\end{figure}

In string diagrams, we will indicate projection maps by drawing a dashed box around a subdiagram, as in Figure \ref{fig:subdiagram}. We stress that the projection maps are only maps of sets, not of spans.

\subsection{Pseudomonoids}\label{sec:pseudo}

Pseudomonoids can be defined in any monoidal bicategory, but because our focus is $\Span$, we use the notation $\times$ for the monoidal product and $\pt$ for the monoidal unit.

\begin{definition}[\cite{DayStreet}, also see \cite{Verdon}]
    
A \emph{pseudomonoid} in a monoidal bicategory consists of
\begin{itemize}
    \item an object $X$,
    \item morphisms $\eta: \pt \to X$ (unit) and $\mu: X \times X \to X$ (multiplication),
    \item invertible $2$-morphisms $a: \mu \circ (\mu \times \id_X) \Rightarrow \mu \circ (\id_X \times \mu)$ (associator), $\ell: \mu \circ (\eta \times \id_X) \Rightarrow \id_X$ (left unitor), and $r: \mu \circ (\id_X \times \eta) \Rightarrow \id_X$ (right unitor),
\end{itemize}
satisfying two coherence conditions, described below.

\end{definition}

We represent $\eta$ and $\mu$, respectively, by the following string diagrams:
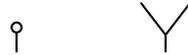
\stringdiagram{
\begin{scope}
    \unit{0}{0}
\end{scope}
\begin{scope}[shift={(6,0)}]
    \multiplication{0}{0}
\end{scope}
}
Note that, in diagrams constructed out of $\eta$ and $\mu$, every string is implicitly labeled by the same object $X$.
The $2$-morphisms $a$, $\ell$ and $r$ can be drawn as morphisms of diagrams:
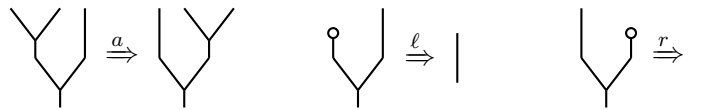
\stringdiagram{
\begin{scope} 
\multiplication{0}{1}
\identity{2}{1}
\multiplication{1}{-1}
\nattrans{3.5}{0}{a}
\identity{5}{1}
\multiplication{7}{1}
\multiplication{6}{-1}
\end{scope}

\begin{scope}[shift={(12,0)}]
\unit{0}{1}
\identity{2}{1}
\multiplication{1}{-1}
\nattrans{3.5}{0}{\ell}
\identity{5}{0}
\end{scope}

\begin{scope}[shift={(22,0)}]
\unit{2}{1}
\identity{0}{1}
\multiplication{1}{-1}
\nattrans{3.5}{0}{r}
\identity{5}{0}
\end{scope}
}

We can now describe the coherence conditions as equations between compositions of $2$-morphisms. The first coherence condition is the \emph{triangle equation}:

	\stringdiagram
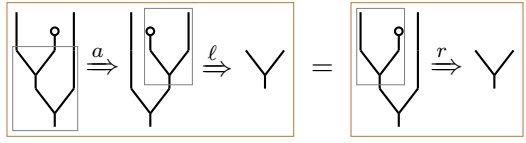
{
        \begin{scope}[scale=0.77]
		\begin{scope}
			\begin{scope}
				\identity{0}{3}
				\unit{2}{3}
				\identity{3}{3}
				\multiplication{1}{1}
				\identity{3}{1}
				\multiplication{2}{-1}
				\highlight{-.2}{-2.2}{3.2}{2.2}
				\nattrans{4.5}{1}{a}
			\end{scope}
			\begin{scope}[shift={(6,0)}]
				\identity{0}{3}
				\unit{1}{3}
				\identity{3}{3}
				\multiplication{2}{1}
				\identity{0}{1}
				\multiplication{1}{-1}
				\highlight{.7}{.2}{3.2}{4.2}
				\nattrans{4.5}{1}{\ell}
			\end{scope}
			\begin{scope}[shift={(12,0)}]
				\multiplication{1}{1}
			\end{scope}
			\enclose{-.5}{-2.5}{14.5}{4.5}
		\end{scope}
		\equals{16}{1}
		\begin{scope}[shift={(18,0)}]
			\begin{scope}
				\identity{0}{3}
				\unit{2}{3}
				\identity{3}{3}
				\multiplication{1}{1}
				\identity{3}{1}
				\multiplication{2}{-1}
				\highlight{-.2}{.2}{2.3}{4.2}
				\nattrans{4.5}{1}{r}
			\end{scope}
			\begin{scope}[shift={(6,0)}]
				\multiplication{1}{1}
			\end{scope}
			\enclose{-.5}{-2.5}{8.5}{4.5}
		\end{scope}
        \end{scope}
	}
The second coherence condition is the \emph{pentagon equation}:
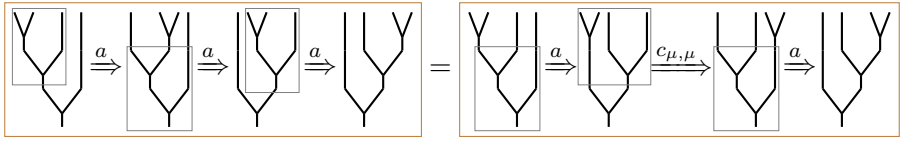
\stringdiagram{
\begin{scope}[scale=0.77]
\begin{scope}
    \begin{scope}
        \slimmultiplication{0}{5}
        \identity{2}{5}
        \identity{3}{5}
        \multiplication{1}{3}
        \identity{3}{3}
        \multiplication{2}{1}
        \nattrans{4.3}{3}{a}
        \highlight{-0.6}{2.2}{2.2}{6.2}
    \end{scope}
    \begin{scope}[shift={(5.6,0)}]
        \slimmultiplication{2}{5}
        \identity{0}{5}
        \identity{3}{5}
        \multiplication{1}{3}
        \identity{3}{3}
        \multiplication{2}{1}
        \nattrans{4.3}{3}{a}
        \highlight{-0.2}{-0.2}{3.2}{4.2}
    \end{scope}
    \begin{scope}[shift={(11.2,0)}]
        \slimmultiplication{1}{5}
        \identity{0}{5}
        \identity{3}{5}
        \multiplication{2}{3}
        \identity{0}{3}
        \multiplication{1}{1}
        \nattrans{4.3}{3}{a}
        \highlight{0.4}{1.8}{3.2}{6.2}
    \end{scope}
    \begin{scope}[shift={(16.8,0)}]
        \slimmultiplication{3}{5}
        \identity{0}{5}
        \identity{1}{5}
        \multiplication{2}{3}
        \identity{0}{3}
        \multiplication{1}{1}
    \end{scope}  
    \enclose{-1}{-.5}{20.8}{6.5}
\end{scope}
\equals{21.8}{3}
\begin{scope}[shift={(23.8,0)}]
    \begin{scope}
        \slimmultiplication{0}{5}
        \identity{2}{5}
        \identity{3}{5}
        \multiplication{1}{3}
        \identity{3}{3}
        \multiplication{2}{1}
        \nattrans{4.3}{3}{a}
        \highlight{-0.2}{-0.2}{3.2}{4.2}
    \end{scope}
    \begin{scope}[shift={(5.8,0)}]
        \slimmultiplication{0}{5}
        \identity{1}{5}
        \identity{3}{5}
        \multiplication{2}{3}
        \identity{0}{3}
        \multiplication{1}{1}
        \nattrans{4.8}{3}{c_{\mu,\mu}}
        \highlight{-0.6}{2.2}{3.2}{6.2}
    \end{scope}
    \begin{scope}[shift={(12.5,0)}]
        \slimmultiplication{3}{5}
        \identity{0}{5}
        \identity{2}{5}
        \multiplication{1}{3}
        \identity{3}{3}
        \multiplication{2}{1}
        \nattrans{4.3}{3}{a}
        \highlight{-0.2}{-0.2}{3.2}{4.2}
    \end{scope}
    \begin{scope}[shift={(18,0)}]
        \slimmultiplication{3}{5}
        \identity{0}{5}
        \identity{1}{5}
        \multiplication{2}{3}
        \identity{0}{3}
        \multiplication{1}{1}
    \end{scope}  
    \enclose{-1}{-.5}{22}{6.5}
\end{scope}
\end{scope}
}
In the above equations, the gray boxes indicate the portion of the diagram where the $2$-morphism is applied. We note that the slide move \eqref{eqn:slide} appears in the pentagon equation. 

\section{\texorpdfstring{$2$-Segal sets correspond to pseudomonoids in $\Span$}{2-Segal sets correspond to pseudomonoids in Span}}
\label{sec:correspondence}

The goal of this section is to prove the following, which is a special case of Stern's correspondence \cite{Stern:span}.
\begin{theorem}\label{thm:correspondence}
    There is a one-to-one correspondence between isomorphism classes of $2$-Segal sets and isomorphism classes of pseudomonoids in $\Span$. 
\end{theorem}
Here, an isomorphism of pseudomonoids is defined as a strong monoidal equivalence. In $\Span$, this notion is rather strict. Specifically, if $X$ and $Y$ are pseudomonoids in $\Span$, an isomorphism between them consists of an isomorphism of sets $X \cong Y$ and isomorphisms of their respective unit and multiplication spans, commuting with the associators and unitors.

In Section \ref{sec:2Segaltopseudomonoid}, we follow the approach in \cite{CMS}*{Section 3.3} to show that one can obtain a pseudomonoid in $\Span$ from a $2$-Segal set. In Section \ref{sec:pseudomonoidto2Segal}, we describe the process of obtaining a $2$-Segal set from a pseudomonoid in $\Span$.

\subsection{\texorpdfstring{From $2$-Segal sets to pseudomonoids in $\Span$}{From 2-Segal sets to pseudomonoids in Span}}\label{sec:2Segaltopseudomonoid}

Let $X_{\bullet}$ be a $2$-Segal set. We can then define a pseudomonoid in $\Span$ as follows:
\begin{itemize}
    \item the underlying object is $X_1$,
    \item the unit and multiplication morphisms are
    \begin{align}\label{eqn:unitandmultiplication}
        \eta&: \pt \leftarrow X_0 \xrightarrow{s_0} X_1, & \mu&: X_1 \times X_1 \xleftarrow{(d_2,d_0)} X_2 \xrightarrow{d_1} X_1.
    \end{align}
\end{itemize}
We still need to define an associator $a$ and unitors $\ell,r$ satisfying the pentagon and triangle equations (see Section \ref{sec:pseudo}). We will see that these naturally arise from a richer structure coming from the $2$-Segal property.

For each $n \geq 0$, we consider the spans
\begin{equation}\label{eqn:nfold}
\mu_n: (X_1)^n \xleftarrow{(e_1,\dots,e_n)} X_n \xrightarrow{e_{\out}} X_1,
\end{equation}
which can be interpreted as $n$-fold product maps. Note that $\mu_0 = \eta$, $\mu_1 = \id_{X_1}$, and $\mu_2 = \mu$.

We now have (at least when $n \geq 2$) two different ways to graphically represent the set $X_n$. In the graphical calculus of Section \ref{sec:gcalc}, $X_n$ is represented by an $(n+1)$-gon. On the other hand, viewing $X_n$ as the apex of the span $\mu_n$, we can represent it by a string diagram with $n$ incoming strings and one outgoing string (where every string is implicitly labeled by $X_1$). There is a simple relationship between these two depictions; the string diagram for $\mu_n$ can be obtained as the dual graph of the $(n+1)$-gon, and vice versa (see Figure \ref{fig:dualgraph})).

\begin{figure}[t]
\begin{center}
\begin{tikzpicture}[scale=0.5]
    \begin{scope} [shift={(0,0)}]
        \draw (-150:1) node[vertex] (0) {}
        -- (-30:1) node[vertex] (2) {} 
        -- (90:1) node[vertex] (1) {} 
        -- cycle;
        \begin{scope}[color=Gray]
        \node[vertex,draw,circle,fill] (A) at (0:0) {};
        \node[vertex] (out) at (270:1.5) {};
        \node[vertex] (1) at (150:1.5) {};
        \node[vertex] (2) at (30:1.5) {};
        \draw (A) -- (out);
        \draw (A) -- (1);
        \draw (A) -- (2);
        \end{scope}
    \end{scope}
    \begin{scope} [shift={(5,0)}]
        \draw (-135:1) node[vertex] (0) {}
        -- (-45:1) node[vertex] (3) {} 
        -- (45:1) node[vertex] (2) {} 
        -- (135:1) node[vertex] (1) {} 
        -- cycle;
        \begin{scope}[color=Gray]
        \node[vertex,draw,circle,fill] (A) at (0:0) {};
        \node[vertex] (out) at (270:1.5) {};
        \node[vertex] (1) at (180:1.5) {};
        \node[vertex] (3) at (0:1.5) {};
        \node[vertex] (2) at (90:1.5) {};
        \draw (A) -- (out);
        \draw (A) -- (1);
        \draw (A) -- (2);
        \draw (A) -- (3);
        \end{scope}
    \end{scope}
    \begin{scope}[shift={(10,0)}]
        \draw (-126:1) node[vertex] (0) {} 
        -- (-54: 1) node[vertex] (4) {} 
        -- (18:1) node[vertex] (3) {} 
        -- (90:1) node[vertex] (2) {} 
        -- (162:1) node[vertex] (1) {} 
        -- cycle;
        \begin{scope}[color=Gray]
        \node[vertex,draw,circle,fill] (A) at (0:0) {};
        \node[vertex] (out) at (270:1.5) {};
        \node[vertex] (1) at (198:1.5) {};
        \node[vertex] (2) at (126:1.5) {};
        \node[vertex] (3) at (54:1.5) {};
        \node[vertex] (4) at (-18:1.5) {};
        \draw (A) -- (out);
        \draw (A) -- (1);
        \draw (A) -- (2);
        \draw (A) -- (3);
        \draw (A) -- (4);        
        \end{scope}
    \end{scope}
\end{tikzpicture}
\end{center}
\caption{The dual graph of the $(n+1)$-gon can be interpreted as the string diagram for the $n$-fold multiplication span \eqref{eqn:nfold}.}
\label{fig:dualgraph}
\end{figure}

This dual graph relationship extends to subdivided polygons. In Section \ref{sec:gcalc}, a subdivision of an $(n+1)$-gon into a $(k+1)$-gon and an $(n-k+2)$-gon is used to represent the pullback of $X_k$ and $X_{n-k+1}$ over the shared edge in $X_1$. Taking the dual graph, we obtain a string diagram that represents a composition of $\mu_k$ and $\mu_{n-k+1}$, defined by exactly the same pullback of $X_k$ and $X_{n-k+1}$; see Figure \ref{fig:dualgraphsubdivide} for an example of this.
\begin{figure}[t]
    \begin{center}
    \begin{tikzpicture}[scale=0.5]
            \begin{scope}
        \draw (240:1) node[vertex] (0) {} 
        -- (300:1) node[vertex] (5) {}
        -- (0:1) node[vertex] (4) {} 
        -- (60:1) node[vertex] (3) {}
        -- (120:1) node[vertex] (2) {}
        -- (180:1) node[vertex] (1) {}    
        -- cycle;
        \draw (1) -- (3);
            \begin{scope}[color=Gray]
        \node[vertex,draw,circle,fill] (A) at (-60:.1) {};
        \node[vertex,draw,circle,fill] (B) at (120:.7) {};        
        \node[vertex] (out) at (270:1.5) {};
        \node[vertex] (1) at (210:1.5) {};
        \node[vertex] (2) at (150:1.5) {};
        \node[vertex] (3) at (90:1.5) {};
        \node[vertex] (4) at (30:1.5) {};
        \node[vertex] (5) at (-30:1.5) {};
        \draw (B) -- (A);
        \draw (B) -- (2);
        \draw (B) -- (3);
        \draw (A) -- (out);
        \draw (A) -- (1);
        \draw (A) -- (4);      
        \draw (A) -- (5);   
        \end{scope}
    \end{scope}
    \end{tikzpicture}
    \end{center}
    \caption{The subdivided hexagon and the string diagram obtained as its dual graph are two different ways of representing the same pullback $X_2 \times_{X_1} X_4$.}
    \label{fig:dualgraphsubdivide}
\end{figure}

This correspondence actually gives a slight enhancement to our interpretation of subdivided polygons, allowing us to consider them as spans, rather than merely as sets. To be more specific, any subdivided $(n+1)$-gon can be viewed as representing a span from $(X_1)^n$ (corresponding to the vertebral edges) to $X_1$ (corresponding to the long edge). A key observation is that the maps \eqref{eqn:subdivide3}, which generate the $2$-Segal functors (see Section \ref{sec:gcalc}), are maps of spans. Thus we can view the $2$-Segal functors as going from the poset category of subdivided $(n+1)$-gons to $\Hom_{\Span}((X_1)^n, X_1)$.
\begin{figure}[t]
\begin{center}
\begin{tikzcd}
&  & (X_1)^3 &  &\\
X_2 \times_{X_1} X_2 \arrow[rru] \arrow[rrd] &  & X_3 \arrow[d, "e_\out"] \arrow[u, "{(e_1,e_2,e_3)}"'] \arrow[rr, LA, "\mathcal{T}_{02}"] \arrow[ll, LA, "\mathcal{T}_{13}"'] &  & X_2 \times_{X_1} X_2 \arrow[llu] \arrow[lld] \\
&  & X_1 & &
\end{tikzcd}
\end{center}
    \caption{An enhancement of the $2$-Segal functor for $n=3$ (see \eqref{eqn:tacos}). The maps $\mathcal{T}_{13}$ and $\mathcal{T}_{02}$ are viewed as elements of $\Hom_{\Span}((X_1)^3, X_1)$, and they are bold to indicate that they are isomorphisms when $X_\bullet$ is $2$-Segal. An important point is that all of the maps in this diagram are defined completely in terms of the simplicial structure, and the commutativity of the diagram follows from simplicial identities. For example, the map in the upper right from $X_2 \times_{X_1} X_2$ to $(X_1)^3$ is given by $(\xi_0, \xi_2) \mapsto (d_2 \xi_2, d_2 \xi_0, d_0 \xi_0)$.}
    \label{fig:2segalfunctspan}
\end{figure}

The $n=3$ case is shown in Figure \ref{fig:2segalfunctspan}. The spans that appear there, in $\Hom_{\Span}((X_1)^3, X_1)$, can be represented by string diagrams as follows:
\stringdiagram{
\begin{scope}[scale=.77]
\begin{scope}
\multiplication{0}{1}
\identity{2}{1}
\multiplication{1}{-1}
\end{scope}
\leftnattrans{4}{0}{\mathcal{T}_{13}}
\begin{scope}[shift={(7,0)}]
\tripleprod{0}{0}    
\end{scope}
\nattrans{10}{0}{\mathcal{T}_{02}}
\begin{scope}[shift={(12,0)}]
\identity{0}{1}
\multiplication{2}{1}
\multiplication{1}{-1}    
\end{scope}
\end{scope}
}
Note that these string diagrams are the dual graphs of the subdivided squares in Figure \ref{fig:assoc3}.

Since $X_\bullet$ is assumed to be $2$-Segal, the maps $\mathcal{T}_{13}$ and $\mathcal{T}_{02}$ are invertible, and we can define an associator 
\[a = \mathcal{T}_{02} (\mathcal{T}_{13})^{-1}.\]
To see that the pentagon equation holds, we consider the image of the $2$-Segal functor for $n=3$. It's a nice exercise for the reader to draw the diagram of dual graphs corresponding to Figure \ref{fig:assoc4} and observe that the commutativity of the diagram implies the pentagon equation.

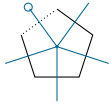
\begin{figure}[t]
    \begin{center}
    \begin{tikzpicture}[scale=0.5]
    \begin{scope}
        \draw (162:1) node[vertex] (1) {} 
        -- (-126:1) node[vertex] (0) {} 
        -- (-54: 1) node[vertex] (4) {} 
        -- (18:1) node[vertex] (3) {} 
        -- (90:1) node[vertex] (2) {};
        \draw[densely dotted] (1) -- (2);
        \begin{scope}[color=Gray]
        \node[vertex,draw,circle,fill] (A) at (0:0) {};
        \node[vertex] (out) at (270:1.5) {};
        \node[vertex] (1) at (198:1.5) {};
        \node[vertex] (2) at (126:1.3) {};
        \node[vertex] (3) at (54:1.5) {};
        \node[vertex] (4) at (-18:1.5) {};
        \draw (A) -- (out);
        \draw (A) -- (1);
        \draw (A) -- (2) circle (3pt);
        \draw (A) -- (3);
        \draw (A) -- (4);        
        \end{scope}
    \end{scope}
    \end{tikzpicture}
    \end{center}
    \caption{The pentagon with a dotted edge and the string diagram obtained as its dual graph are two different ways of representing $X_0 \times_{X_1} X_4$, which is canonically isomorphic to $\{\omega \in X_4 \suchthat e_2 \omega \in s_0(X_0)\}$.}
    \label{fig:dualdotted}
\end{figure}

To define unitors, we further extend the dual graph correspondence to allow for subdivided $(n+1)$-gons where one or more of the edges corresponding to the vertebra $e_1, \dots, e_n$ are dotted. Recall that a dotted edge is used to impose the restriction that a $1$-simplex is degenerate. In the dual graph correspondence, we make the rule that the dual of a dotted edge is the string diagram for $\eta$. The apex of the span represented by the resulting string diagram is canonically isomorphic to the set represented by the partially-dotted subdivided polygon; see Figure \ref{fig:dualdotted}.

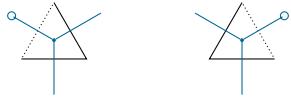
\begin{figure}[t]
\begin{center}
\begin{tikzpicture}[scale=0.5]
    \begin{scope} [shift={(0,0)}]
        \draw (-150:1) node[vertex] (0) {}
        -- (-30:1) node[vertex] (2) {} 
        -- (90:1) node[vertex] (1) {};
        \draw[densely dotted] (0) -- (1);
        \begin{scope}[color=Gray]
        \node[vertex,draw,circle,fill] (A) at (0:0) {};
        \node[vertex] (out) at (270:1.5) {};
        \node[vertex] (1) at (150:1.3) {};
        \node[vertex] (2) at (30:1.5) {};
        \draw (A) -- (out);
        \draw (A) -- (1) circle (3pt);
        \draw (A) -- (2);
        \end{scope}
    \end{scope}
    \begin{scope} [shift={(5,0)}]
        \draw (-30:1) node[vertex] (2) {}
        --(-150:1) node[vertex] (0) {}
        -- (90:1) node[vertex] (1) {};
        \draw[densely dotted] (1) -- (2);
        \begin{scope}[color=Gray]
        \node[vertex,draw,circle,fill] (A) at (0:0) {};
        \node[vertex] (out) at (270:1.5) {};
        \node[vertex] (1) at (150:1.5) {};
        \node[vertex] (2) at (30:1.3) {};
        \draw (A) -- (out);
        \draw (A) -- (1);
        \draw (A) -- (2) circle (3pt);
        \end{scope}
    \end{scope}

\end{tikzpicture}
\end{center}
\caption{The partially-dotted triangle on the left represents $\{\omega \in X_2 \suchthat e_1 \omega \in s_0(X_0)\}$. The isomorphism with $s_0(X_1) \cong X_1$ defines the left unitor $\ell$. Similarly, the right unitor $r$ is obtained via the triangle on the right.
}
\label{fig:2degen}
\end{figure}

\begin{figure}[t]
\begin{center}
\begin{tikzpicture}[scale=0.5]

    \begin{scope} [shift={(-5,0)}]
        \draw (135:1) node[vertex] (1) {} 
        -- (-135:1) node[vertex] (0) {}
        -- (-45:1) node[vertex] (3) {} 
        -- (45:1) node[vertex] (2) {};
        \draw[densely dotted] (1) -- (2);
        \draw (0) -- (2);
        \begin{scope}[color=Gray]
        \node[vertex,draw,circle,fill] (A) at (135:.4) {};
        \node[vertex,draw,circle,fill] (B) at (-45:.4) {};
        \node[vertex] (out) at (270:1.5) {};
        \node[vertex] (1) at (180:1.5) {};
        \node[vertex] (3) at (0:1.5) {};
        \node[vertex] (2) at (90:1.3) {};
        \draw (A) -- (B);
        \draw (A) -- (1);
        \draw (A) -- (2) circle (3pt);
        \draw (B) -- (3);
        \draw (B) -- (out);
        \end{scope}
    \end{scope}

    \begin{scope} [shift={(0,0)}]
        \draw (135:1) node[vertex] (1) {} 
        -- (-135:1) node[vertex] (0) {}
        -- (-45:1) node[vertex] (3) {} 
        -- (45:1) node[vertex] (2) {};
        \draw[densely dotted] (1) -- (2);
        \begin{scope}[color=Gray]
        \node[vertex,draw,circle,fill] (A) at (0:0) {};
        \node[vertex] (out) at (270:1.5) {};
        \node[vertex] (1) at (180:1.5) {};
        \node[vertex] (3) at (0:1.5) {};
        \node[vertex] (2) at (90:1.3) {};
        \draw (A) -- (out);
        \draw (A) -- (1);
        \draw (A) -- (2) circle (3pt);
        \draw (A) -- (3);
        \end{scope}
    \end{scope}
    \begin{scope} [shift={(5,0)}]
        \draw (135:1) node[vertex] (1) {} 
        -- (-135:1) node[vertex] (0) {}
        -- (-45:1) node[vertex] (3) {} 
        -- (45:1) node[vertex] (2) {};
        \draw[densely dotted] (1) -- (2);
        \draw (1) -- (3);
        \begin{scope}[color=Gray]
        \node[vertex,draw,circle,fill] (A) at (45:.4) {};
        \node[vertex,draw,circle,fill] (B) at (-135:.4) {};
        \node[vertex] (out) at (270:1.5) {};
        \node[vertex] (1) at (180:1.5) {};
        \node[vertex] (3) at (0:1.5) {};
        \node[vertex] (2) at (90:1.3) {};
        \draw (A) -- (B);
        \draw (A) -- (3);
        \draw (A) -- (2) circle (3pt);
        \draw (B) -- (1);
        \draw (B) -- (out);
        \end{scope}
    \end{scope}
    \begin{scope} [shift={(0,-5)}]
        \draw (-150:1) node[vertex] (0) {}
        -- (-30:1) node[vertex] (2) {} 
        -- (90:1) node[vertex] (1) {} 
        -- cycle;
        \begin{scope}[color=Gray]
        \node[vertex,draw,circle,fill] (A) at (0:0) {};
        \node[vertex] (out) at (270:1.5) {};
        \node[vertex] (1) at (150:1.5) {};
        \node[vertex] (2) at (30:1.5) {};
        \draw (A) -- (out);
        \draw (A) -- (1);
        \draw (A) -- (2);
        \end{scope}
    \end{scope}
\nattrans{2.5}{0}{\mathcal{T}_{02}}
\leftnattrans{-2.5}{0}{\mathcal{T}_{13}}

\node[rotate=90] at (0,-2.5) {$\Longrightarrow$};
\node[anchor=west] at (0,-2.5) {$s_1$};
\end{tikzpicture}
\end{center}
\caption{The maps that appear in the triangle equation are $a = \mathcal{T}_{02} \circ \mathcal{T}_{13}^{-1}$, $\ell = s_1^{-1} \circ \mathcal{T}_{02}^{-1}$, $r = s_1^{-1} \circ \mathcal{T}_{13}^{-1}$. The triangle equation $r = \ell \circ a$ follows.}
\label{fig:triangleidentity}
\end{figure}
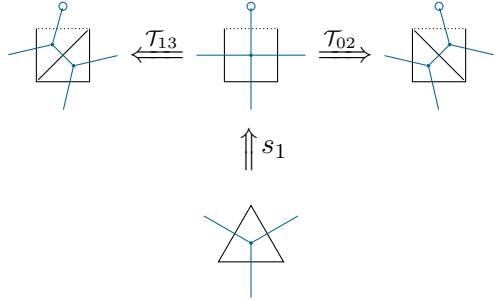

In Section \ref{sec:gcalc}, it was observed that, as a result of the unitality property of Remark \ref{remark:unitality}, $\omega \in X_n$ satisfies $e_{i+1} \omega \in s_0(X_0)$ if and only if $\omega \in s_i(X_{n-1})$. As a result, the sets represented by the partially-dotted triangles in Figure \ref{fig:2degen} are canonically isomorphic to $s_0(X_1)$ and $s_1(X_1)$. We can then define unitors $\ell$ and $r$ as the isomorphisms $s_0(X_1) \isoto X_1$ and $s_1(X_1) \isoto X_1$ obtained as restrictions of $d_1^2$. The triangle equation is then an immediate consequence of the diagram in Figure \ref{fig:triangleidentity}.

\subsection{From pseudomonoids in Span to 2-Segal sets}\label{sec:pseudomonoidto2Segal}

Let $X$ be a pseudomonoid in $\Span$ with unit morphism $\eta: \pt \to X$ and multiplication morphism $\mu: X \times X \to X$. By repeatedly using $\mu$, we define the (left-to-right) $n$-fold multiplication morphism $\mu_n: X^n \to X$, given by
\[ \mu \circ (\mu \times \id_X) \circ (\mu \times \id_{X^2}) \circ \cdots \circ (\mu \times \id_{X^{n-2}})\]
for $n \geq 2$, with $\mu_1 = \id_X$ and $\mu_0 = \eta$. 

\begin{figure}[t]
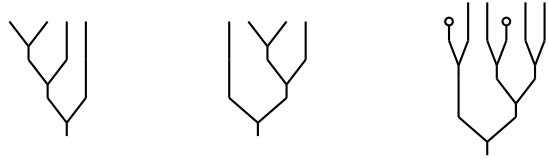

\stringdiagram{
\begin{scope}[scale=.77]
\begin{scope}
\multiplication{0}{2}
\identity{2}{2}
\identity{3}{2}
\multiplication{1}{0}
\identity{3}{0}
\multiplication{2}{-2}
\end{scope}
\begin{scope}[shift={(10,0)}]
\identity{.5}{2}
\multiplication{2.5}{2}
\identity{4.5}{2}
\identity{.5}{0}
\multiplication{3.5}{0}
\widemultiplication{2}{-2}
\end{scope}
\begin{scope}[shift={(22,-1)}]
\unit{0}{4}
\identity{1}{4}
\identity{2}{4}
\unit{3}{4}
\identity{4}{4}
\identity{5}{4}
\slimmultiplication{.5}{2}
\slimmultiplication{2.5}{2}
\slimmultiplication{4.5}{2}
\identity{.5}{0}
\multiplication{3.5}{0}
\widemultiplication{2}{-2}
\end{scope}
\end{scope}
}
    \caption{On the left, the string diagram for $\mu_4$. In the middle and right, two other trees with four open leaves, representing morphisms that are canonically $2$-isomorphic to $\mu_4$.}
    \label{fig:trees}
\end{figure}

More generally, the morphisms that can be constructed out of $\eta$ and $\mu$ are represented by string diagrams whose connected components are planar binary rooted trees with leaves that are either open or ``capped'' by $\eta$. An important fact is the coherence theorem for pseudomonoids \cites{lack:coherent,Verdon}, which implies that any morphism represented by a tree with $n$ open leaves is $2$-isomorphic to $\mu_n$, via a unique $2$-isomorphism constructed out of the associator, unitors, and the structure $2$-morphisms for the underlying monoidal bicategory; see Figure \ref{fig:trees}.

\begin{figure}[t]
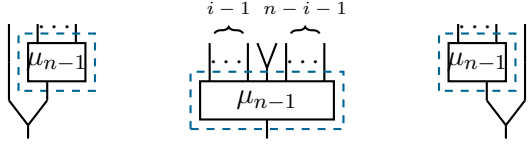

       \stringdiagram{
       \begin{scope}[scale=.77]
        \begin{scope}
        \identity{0}{3}
        \identity{0}{5}
        \halfidentity{1.5}{6}
        \strDots{2.5}{5.8}
        \halfidentity{3.5}{6}
        \boxA{1}{3}{\mu_{n-1}}
        \dbox{0.5}{2.5}
        \halfidentity{2}{3}
        \multiplication{1}{1}
        \end{scope}

        \begin{scope}[shift={(10,1)}]
            \halfidentity{3.5}{0}
            \identity{0.5}{3}
            \identity{2.5}{3}
            \strDots{1.5}{3}
            \slimmultiplication{3.5}{3}
            \identity{4.5}{3}
            \strDots{5.5}{3}
            \identity{6.5}{3}
            \longboxA{0}{0}{\mu_{n-1}}
            \longdbox{-0.5}{-0.5}
            \strOverBrace{0.75}{4.5}{\tiny $i-1$}
            \strOverBrace{4.75}{4.5}{\tiny $n-i-1$}
        \end{scope}

        \begin{scope}[shift={(27,0)}]
        \identity{0}{3}
        \identity{0}{5}
        \halfidentity{-1.5}{6}
        \strDots{-2.5}{5.8}
        \halfidentity{-3.5}{6}
        \boxA{-4}{3}{\mu_{n-1}}
        \dbox{-4.5}{2.5}
        \halfidentity{-2}{3}
        \multiplication{-1}{1}
        \end{scope}
        \end{scope}
        }
    \caption{String diagrams depicting the face maps $d_0^n$ (left), $d_i^n$ for $0<i<n$ (middle), and $d_n^n$ (right). The face maps are defined by using the canonical $2$-isomorphism from $\mu_n$ to the morphism shown here, then projecting to the $\mu_{n-1}$ component.}
    \label{fig:facemapgraphic}
\end{figure}

We will now construct a simplicial set $X_\bullet$ associated to the pseudomonoid $X$. For each $n$, let $X_n$ be the apex of the span $\mu_n$. The inner face maps $d_i^n: X_n \to X_{n-1}$, $0 < i < n$, are defined by using the canonical $2$-isomorphism
\[ \mu_n \Rightarrow \mu_{n-1} \circ (\id_{X^{i-1}} \times \mu \times \id_{X^{n-i-1}})\]
and then projecting onto the $\mu_{n-1}$ component (recall that projection maps were described in Section \ref{sec:projection}). The outer face maps $d_0^n$ and $d_n^n$ are defined by using the canonical $2$-isomorphisms
\begin{align*} 
\mu_n &\Rightarrow \mu \circ (\id_X \times \mu_{n-1}), & \mu_n &\Rightarrow \mu \circ (\mu_{n-1} \times \id_X),
\end{align*}
and then projecting onto the $\mu_{n-1}$ component. The string diagrams depicting the face maps are shown in Figure \ref{fig:facemapgraphic}.

\begin{figure}[t]
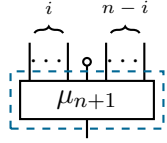

       \stringdiagram{
       \begin{scope}[scale=.77]
            \begin{scope}
            \halfidentity{3.5}{0}
            \identity{0.5}{3}
            \identity{2.5}{3}
            \strDots{1.5}{3}
            \unit{3.5}{3}
            \identity{4.5}{3}
            \strDots{5.5}{3}
            \identity{6.5}{3}
            \longboxA{0}{0}{\mu_{n+1}}
            \longdbox{-0.5}{-0.5}
            \strOverBrace{0.75}{4.5}{\tiny $i$}
            \strOverBrace{4.75}{4.5}{\tiny $n-i$}
            \end{scope}
        \end{scope}
        }
    \caption{String diagram depicting the degeneracy maps $s_i^n$. The degeneracy maps are defined by using the canonical $2$-isomorphism from $\mu_n$ to the morphism shown here, then projecting to the $\mu_{n+1}$ component.}
    \label{fig:degenmapgraphic}
\end{figure}

The degeneracy maps $s_i^n: X_n \to X_{n+1}$ are defined by using the canonical $2$-isomorphism
\[ \mu_n \Rightarrow \mu_{n+1} \circ (\id_{X^i} \times \mu \times \id_{X^{n-i}})\]
and then projecting onto the $\mu_{n+1}$ component; see Figure \ref{fig:degenmapgraphic}.

At this point, one could directly use the above definitions to prove that the simplicial identities hold. Instead, our approach will be to recognize that planar binary rooted trees (with some leaves ``capped'') are exactly the dual graphs of triangulated polygons (with some vertebral edges dotted). Under this correspondence, the above definitions of the face and degeneracy maps correspond precisely to the graphical calculus described in Section \ref{sec:gcalc}. One can then deduce the simplicial identities from the graphical calculus, as in Figure \ref{fig:faceface}.

The $2$-Segal conditions for $X_\bullet$ are almost immediate from this construction, since the string diagram representing $\mu_n$ is the dual graph of the standard left-triangulated $(n+1)$-gon (consisting of all the edges connecting vertex $0$ to the other vertices), and the coherence theorem provides canonical isomorphisms to the sets represented by any other triangulation.

\section{\texorpdfstring{From $2$-Segal sets to associative algebras}{From 2-Segal sets to associative algebras}}
\label{sec:hallalgebra}

In this section, we review the linearization functor taking spans to linear maps of vector spaces, and we describe various constructions that make use of the linearization functor to produce an associative (co)algebra from a $2$-Segal set that satisfies finiteness conditions. These constructions appear in \cites{cooper-young, Dyckerhoff-Kapranov:Higher, GKT1} in the more general context of $2$-Segal spaces (also see \cite{bergner2024combinatorial}). From the perspective of \cite{BaezGroupoid}, where the linearization functor is defined in the more general setting of spans of groupoids, the (co)algebras produced via these constructions could be considered \emph{decategorifications} of the $2$-Segal set.

At the end of the section, we give several examples of algebras that arise from $2$-Segal sets and that may be of interest to readers of different backgrounds and motivations.

Throughout this section, $\kk$ is an arbitrary field.

\subsection{The linearization functor}\label{sec:linearization}

The process of obtaining linear maps from spans was briefly sketched in Section \ref{sec:introduction}, but we provide additional details here.

Given a set $X$, we can form the vector space $\kk[X]$ generated by elements of $X$. Any map of sets $f:X \to Y$ naturally extends to a linear map $f_*: \kk[X] \to \kk[Y]$. If $f$ has finite preimages, then there is also the map $f^*: \kk[Y] \to \kk[X]$, where
\[ f^*(y) = \sum_{x \in f^{-1}(y)} x.\]
Given a span $f: X \fromlabel{f_1} A \tolabel{f_2} Y$ where $f_1$ has finite preimages, we can obtain the linear map 
\[f_\diamond = (f_2)_* (f_1)^*: \kk[X] \to \kk[Y].\] 
The following proposition (which we leave as an exercise for the reader) lists the key properties that make this assignment functorial.
\begin{proposition}\label{prop:linearization}
\begin{enumerate}
    \item  If $f: X \fromlabel{f_1} A \tolabel{f_2} Y$ and $f': X \fromlabel{f_1'} A' \tolabel{f_2'} Y$ are isomorphic spans (i.e.\ if there exists an isomorphism $A \isoto A'$ that commutes with the maps to $X$ and $Y$) such that $f_1$ and $f_1'$ have finite preimages, then $f_\diamond = f'_\diamond$.
    \item Let
    \[\begin{tikzcd}
        C \arrow[d, "\tilde{h}_1"'] \arrow[r, "\tilde{g}_2"] & B \arrow[d, "h_1"] \\
A \arrow[r, "g_2"] & Y 
    \end{tikzcd}\]
    be a pullback square such that $h_1$ has finite preimages. Then $\tilde{h}_1$ has finite preimages, and
    \[ (h_1)^* (g_2)_* = (\tilde{g}_2)_* (\tilde{h}_1)^*.\]
\end{enumerate}
\end{proposition}

Recall that $\Span_1$ is the homotopy category of the bicategory $\Span$. More explicitly, its objects are sets, and its morphisms are isomorphism classes of spans. We use $\Span_1^\tame$ to denote the subcategory of $\Span_1$ consisting of spans where the left leg has finite preimages. As a direct consequence of Proposition \ref{prop:linearization}, we have the following.
\begin{corollary}\label{cor:linearization}
The assignment $X \mapsto \kk[X]$, $f \mapsto f_\diamond$ gives a functor from $\Span_1^\tame$ to $\vect_\kk$.    
\end{corollary}
We refer to the functor in Corollary \ref{cor:linearization} as the \emph{linearization functor}. We leave the reader to verify that it is in fact a monoidal functor from $(\Span_1^\tame, \times)$ to $(\vect_\kk, \otimes)$.

\subsection{Hall algebras}\label{sec:hall}

Let $X_\bullet$ be a $2$-Segal set such that
\begin{enumerate}
    \item $X_0$ is finite,
    \item $T_2 = (d_2,d_0):X_2 \to X_1 \times X_1$ has finite preimages.
\end{enumerate}

\begin{definition}\label{def:hall}
The \emph{Hall algebra} of $X_\bullet$ is $\kk[X_1]$, with multiplication given by
\[ x \cdot y = \sum_{\omega \in T_2^{-1}(x,y)} d_1 \omega\]
for $x,y \in X_1$. The multiplicative identity element is given by
\[ 1 = \sum_{u \in X_0} s_0 u.\]
\end{definition}

The construction of the Hall algebra can be explained, making use of Theorem \ref{thm:correspondence}, as follows. The pseudomonoid in $\Span$ corresponding to $X_\bullet$ induces a monoid object in $\Span_1$. The multiplication and unit morphisms for the monoid are given by the spans \eqref{eqn:unitandmultiplication}. The conditions given above, that $X_0$ is finite and $T_2$ has finite preimages, ensure that the multiplication and unit are in $\Span_1^\tame$. Applying the linearization functor, we obtain a monoid in $\vect_\kk$ (i.e.\ an associative and unital algebra), and this is the Hall algebra.

In addition to explaining the formulas in Definition \ref{def:hall}, this conceptual explanation makes the associativity and unitality of the Hall algebra immediate.

\subsection{Incidence coalgebras and algebras}\label{sec:incidence}

Let $X_\bullet$ be a $2$-Segal set such that $d_1^2$ has finite preimages. 

\begin{definition}
The \emph{incidence coalgebra} of $X_\bullet$ is $\kk[X_1]$, with comultiplication $\Delta: \kk[X_1] \to \kk[X_1] \otimes \kk[X_1]$ given by
\[ \Delta(x) = \sum_{\omega \in d_1^{-1}(\omega)} d_2 \omega \otimes d_0 \omega\]
for $x \in X_1$. The counit $\varepsilon : \kk[X_1] \to \kk$ is given by
\[ \varepsilon(x) = \begin{cases}
1 & x \in s_0(X_0), \\
0 & x \notin s_0(X_0).
\end{cases}
\]
\end{definition}

The construction of the incidence coalgebra is explained as follows. The bicategory $\Span$ has a dagger structure that reverses spans, turning a span $X \leftarrow A \rightarrow Y$ into $Y \leftarrow A \rightarrow X$. Applying the dagger to a pseudomonoid, we obtain a \emph{pseudocomonoid}\footnote{A pseudocomonoid is defined similarly to a pseudomonoid, but with the morphisms going in the opposite direction. Thus the string diagrams in the definition of pseudocomonoid are reflections (over a horizontal line) of the diagrams in Section \ref{sec:pseudo}.} in $\Span$, which induces a comonoid object in $\Span_1$. The comultiplication and counit morphisms for the comonoid are given by the reverses of the spans \ref{eqn:unitandmultiplication}. The counit is always in $\Span_1^\ell$, and the requirement that $d_1^2$ have finite preimages ensures that the comultiplication is also in $\Span_1^\ell$. Applying the linearization functor, we obtain a comonoid in $\vect_\kk$ (i.e.\ a coalgebra), and this is the incidence coalgebra.

Let $C(X_1) = \Hom(\kk[X_1],\kk)$ be the space of $\kk$-valued functions on $X_1$. As the dual space to a coalgebra, it possesses a natural algebra structure, given by the convolution product
\[
\begin{tikzcd}
{\kk[X_1]} \arrow[r, "\Delta"] \arrow[rrr, "\psi_1*\psi_2", bend right] & {\kk[X_1]\otimes \kk[X_1]} \arrow[r, "\psi_1 \otimes \psi_2"] & \kk \otimes \kk \arrow[r] & \kk
\end{tikzcd}
\]
for $\psi_1, \psi_2 \in C(X_1)$. More explicitly, we have
\begin{equation}\label{eqn:convolution}
   (\psi_1*\psi_2)(x) = \sum_{\omega \in d_1^{-1}(x)} \psi_1(d_2 \omega)\psi_2(d_0 \omega). 
\end{equation} 
The algebra $C(X_1)$ is called the \emph{incidence algebra} of $X_\bullet$.

\begin{remark}\label{remark:hall-incidence}
We emphasize that the Hall algebra and incidence (co)algebra constructions have different finiteness conditions, so it is possible that one is defined when the other is not. When both are defined, the map $\kk[X_1] \to C(X_1)$, taking $x \in X_1$ to its characteristic function, is an injective algebra homomorphism, allowing us to view $\kk[X_1]$ as the subalgebra of $C(X_1)$ consisting of functions with finite support. In particular, when $X_2$ is finite, then $\kk[X_1]$ and $C(X_1)$ are both defined and naturally isomorphic.
\end{remark}

\subsection{Hall and incidence algebras of partial categories}

Recall from Section \ref{sec:hall} that the Hall algebra of a $2$-Segal set $X_\bullet$ is defined when $X_0$ is finite and $T_2$ has finite preimages. If $\cat$ is a partial category, then its nerve $N\cat_\bullet$ automatically satisfies the latter condition (see Proposition \ref{prop:partialcatnerve}), so its Hall algebra is defined whenever $\cat$ has a finite number of objects\footnote{If $\cat$ has infinitely many objects, one can still define an associative product, but not the identity element.}. In this case, the Hall algebra, as a vector space, is generated by morphisms of $\cat$, with multiplication given by
\[ fg = \begin{cases}
    g \circ f & \mbox{ if $g \circ f$ is defined},\\
    0 & \mbox{ otherwise},
\end{cases}
\]
and with identity
\[ 1 = \sum_{x \in \Ob(\cat)} \id_x.\]

On the other hand, the incidence (co)algebra of a $2$-Segal set is defined when $d_1^2$ has finite preimages (see Section \ref{sec:incidence}). In the case of $N\cat_\bullet$ this translates to the condition that, for every morphism $h$, there are a finite number of decompositions of the form $h = g \circ f$. This property (in the case where $\cat$ is a category) is called ``finitely finite'' in \cite{leinster} and ``finite decompositions of degree $2$'' in \cite{leroux}. 

In this case, we denote the incidence algebra as
\[ I(\cat) = \{\psi: \Mor(\cat) \to \kk\}\]
with the product
\[ (\psi_1 * \psi_2)(h) = \sum_{h = g \circ f} \psi_1(f)\psi_2(g).\]

These constructions include many well-known algebras as special cases:
\begin{enumerate}

    \item (Group algebras) Let $G$ be a group, viewed as a category with one object. Then the Hall algebra is the (opposite of the) usual group algebra $\kk[G]$.

    When $G$ is finite, the incidence algebra $I(G)$ is defined, with the well-known convolution product
    \[ (\psi_1*\psi_2)(h) = \sum_{g \in G} \psi_1(g^{-1}h) \psi_2(g).\]
    
    \item (Poset incidence algebras) Let $P$ be a poset, viewed as a category. In order for the Hall algebra to be defined, $P$ must be finite. But the incidence algebra is defined under the weaker condition of \emph{local finiteness}, that, for $x,y \in P$, there exists a finite number of $z$ such that $x \leq z \leq y$.

    In this case, the incidence algebra $I(P)$ is the space of functions on $\{(x,y) \in P\times P \suchthat x \leq y\}$, with the product
    \begin{equation} \label{eqn:posetproduct}
    (\psi_1*\psi_2)(x,y) = \sum_{x \leq z \leq y} \psi_1(x,z) \psi_2(z,y),
    \end{equation}
    which is the product for the incidence algebra of the poset, as defined by Rota \cite{rota:comb1}.

    \item (Dirichlet product) As a particular example of a poset, let $P = \N_+$, the set of positive natural numbers, with the partial ordering given by divisibility. Then the incidence algebra $I(\N_+)$ is the space of functions on $\{(x,y) \in \N_+ \times \N_+ \suchthat x \mbox{ divides } y\}$.

    Given an arithmetic function $\psi: \N_+ \to \kk$, we can define $\hat{\psi} \in I(\N_+)$ by $\hat{\psi}(x,y) = \psi(\frac{y}{x})$. The map $\psi \mapsto \hat{\psi}$ allows us to identify the space of arithmetic functions with a subspace of $I(\N_+)$. Furthermore, one can show that the arithmetic functions are closed under the product \eqref{eqn:posetproduct}, and the restriction coincides with the classical Dirichlet product.
    
    \item For a positive integer $L$, consider the partial monoid $\mathcal{L} = \{1, x, x^2, \dots, x^L\}$ of Example \ref{ex:partialmonoid}. Since $\mathcal{L}$ is finite, the Hall algebra and incidence algebra are both defined and isomorphic to each other. They can be identified with the quotient of the polynomial algebra $\kk[x]$ by the ideal $\langle x^{L+1} \rangle$.
\end{enumerate}

\subsection{More examples}
We end the paper with descriptions of the algebras arising from the examples of $2$-Segal sets in Section \ref{sec:moreexamples}.

\begin{example}
Let $A$ be a set, and let $X_\bullet$ be the $2$-Segal set arising from $A$ in Example \ref{ex:powerset}. Then $X_0 = \{A\}$ is a one-point set, and
\begin{align*}
    X_1 &= \{(B_0,B_1) \in P(A)^2 \suchthat B_0 \cup B_1 = A\},\\
    X_2 &= \{(B_0,B_1,B_2) \in P(A)^3 \suchthat B_0 \cup B_1 \cup B_2 = A\}.
\end{align*}
The face maps from $X_2$ to $X_1$ are given by
\begin{align*}
    d_0^2(B_0,B_1,B_2) &= (B_0 \cup B_1, B_2),\\
    d_1^2(B_0,B_1,B_2) &= (B_0, B_1 \cup B_2),\\
    d_2^2(B_0,B_1,B_2) &= (B_0 \cup B_2, B_1).
\end{align*}

Recall that the Hall algebra construction requires $T_2$ to have finite preimages, and the incidence algebra construction requires $d_1$ to have finite preimages.
As noted in Example \ref{ex:powerset}, $T_2(B,A,A) = ((A,A),(A,A))$ for all $B \in P(A)$. Additionally, $d_1^2(A,A,B) = (A,A)$ for all $B \in P(A)$. From these equations we can see that neither of the finiteness conditions is satisfied when $A$ is infinite. 

Thus we will assume that $A$ is finite, in which case the Hall and incidence algebra are both defined and isomorphic to each other (see Remark \ref{remark:hall-incidence}). If $|A|=k$, then $|X_1| = 3^k$, so the Hall algebra $\kk[X_1]$ is $3^k$-dimensional.

To obtain a nice description of the Hall algebra, it will help to introduce notation for certain elements of $X_1$. For any $B \in P(A)$, write
\begin{align*}
    x_B &= (A,B), & y_B = (B^c, B),
\end{align*}
where $B^c$ is the complement of $B$ in $A$. The Hall algebra products involving these elements satisfy the following rules:
\begin{itemize}
    \item If $B_1 \cap B_2 = \emptyset$, then $x_{B_1} x_{B_2} = x_{B_1 \cup B_2}$. This follows from the fact that $\omega = (A, B_1, B_2)$ is the unique solution to $T_2 \omega = (x_{B_1},x_{B_2})$.
    \item If $B_1 \cap B_2 = \emptyset$, then $y_{B_1} y_{B_2} = y_{B_1 \cup B_2}$. This follows from the fact that $\omega = ((B_1 \cup B_2)^c, B_1, B_2)$ is the unique solution to $T_2 \omega = (y_{B_1},y_{B_2})$.
    \item If $B_1 \cap B_2 = \emptyset$, then $x_{B_1} y_{B_2} = (B_2^c, B_1 \cup B_2)$. This follows from the fact that $\omega = (B_2^c,B_1,B_2)$ is the unique solution to $T_2 \omega = (x_{B_1},y_{B_2})$.
\end{itemize}
If we label the elements of $A$ as $a_1,\dots, a_k$ and write $x_i = x_{\{a_i\}}$, $y_i = y_{\{a_i\}}$, then the above rules imply that
\begin{align*}
    x_B &= \prod_{a_i \in B} x_i, & y_B &= \prod_{a_i \in B} y_i, 
\end{align*}
Every $(C_0,C_1) \in X_1$ can be written as $x_{B_1} y_{B_2} = (B_2^c, B_1 \cup B_2)$ by setting $B_2 = C_0^c$ and $B_1 = C_0 \cap C_1$. Thus, as an algebra, $\kk[X_1]$ is generated by the elements $x_i$ and $y_i$.

We can establish relations satisfied by the generators by making the following observations:
\begin{itemize}
    \item $x_i^2 = x_i + y_i$ for all $i$. This follows from the fact that $\omega = (A, \{a_i\},\{a_i\})$ and $\omega = (\{a_i\}^c, \{a_i\}, \{a_i\})$ are the only solutions to $T_2 \omega = (x_i,x_i)$.
    \item $x_i y_i = 0$ for all $i$. This follows from the fact that there are no solutions to $T_2 \omega = (x_i,y_i)$.
\end{itemize}
From these relations, we have that $y_i = x_i^2 - x_i$ (so the $y_i$ are unnecessary in the generating set) and $x_i^3 - x_i^2 = x_i y_i = 0$. We conclude that there is a surjective algebra morphism
\[ \kk[x_1,\dots,x_k]/\langle x_i^3 - x_i^2 \rangle \to \kk[X_1],\]
and dimension counting shows that this is in fact an isomorphism.
\end{example}

\begin{example}
Let $X_\bullet$ be the $2$-Segal set in Example \ref{ex:exponential}. As described there, we have $X_0 = \{0\}$ and $X_1 \cong \N$. Additionally, we can make the identification
\[ X_2 \cong \{ (m,S) \suchthat S \subseteq \underline{m}\},\]
with face maps given by
\begin{align*}
    d_0(m,S) &= m - |S|, & d_1(m,S) &= m, & d_2(m,S) = |S|.
\end{align*}
The incidence algebra $C(X_1)$ consists of functions $\psi: \N \to \kk$, and \eqref{eqn:convolution} becomes
\begin{align*}
   ( \psi_1 * \psi_2)(m) &= \sum_{S \subseteq \underline{m}} \psi_1(|S|) \psi_2(m-|S|) \\
   &= \sum_{k=0}^m \binom{m}{k} \psi_1(k) \psi_2(m-k).
\end{align*}
This is the Cauchy product arising from exponential generating functions, i.e.\ it is characterized by the equation
\[ \sum_{m=0}^\infty \frac{(\psi_1*\psi_2)(m)}{m!} x^m = \left(\sum_{k=0}^\infty \frac{\psi_1(k)}{k!} x^k \right)\left(\sum_{\ell=0}^\infty \frac{\psi_2(\ell)}{\ell!} x^\ell \right).\]

The Hall algebra of $X_\bullet$ is also defined and can be identified with the subalgebra of $C(X_1)$ consisting of functions with finite support (see Remark \ref{remark:hall-incidence}).
\end{example}

\bibliography{Comm}

\begin{bibdiv}
\begin{biblist}

\bib{BaezDolan:fromfinitesets}{incollection}{
      author={Baez, John~C.},
      author={Dolan, James},
       title={From finite sets to {F}eynman diagrams},
        date={2001},
   booktitle={Mathematics unlimited---2001 and beyond},
   publisher={Springer, Berlin},
       pages={29\ndash 50},
}

\bib{BaezGroupoid}{article}{
      author={Baez, John~C.},
      author={Hoffnung, Alexander~E.},
      author={Walker, Christopher~D.},
       title={Higher dimensional algebra {VII}: {G}roupoidification},
        date={2010},
        ISSN={1201-561X},
     journal={Theory Appl. Categ.},
      volume={24},
       pages={No. 18, 489\ndash 553},
}

\bib{bergner2024combinatorial}{article}{
      author={Bergner, Julia~E.},
       title={Combinatorial examples and applications of 2-{S}egal sets},
        date={2026},
     journal={Bolet\'in de Matem\'aticas},
      volume={32},
      number={2},
       pages={31\ndash 51},
}

\bib{BOORS}{article}{
      author={Bergner, Julia~E.},
      author={Osorno, Ang\'{e}lica~M.},
      author={Ozornova, Viktoriya},
      author={Rovelli, Martina},
      author={Scheimbauer, Claudia~I.},
       title={2-{S}egal sets and the {W}aldhausen construction},
        date={2018},
        ISSN={0166-8641},
     journal={Topology Appl.},
      volume={235},
       pages={445\ndash 484},
         url={https://doi.org/10.1016/j.topol.2017.12.009},
}

\bib{CKWW}{article}{
      author={Carboni, A.},
      author={Kelly, G.~M.},
      author={Walters, R. F.~C.},
      author={Wood, R.~J.},
       title={Cartesian bicategories {II}},
        date={2007},
     journal={Theory Appl. Categ.},
      volume={19},
       pages={93\ndash 124},
}

\bib{leroux}{article}{
      author={Content, Mireille},
      author={Lemay, Fran\c{c}ois},
      author={Leroux, Pierre},
       title={Cat\'egories de {M}\"obius et fonctorialit\'es: un cadre g\'en\'eral pour l'inversion de {M}\"obius},
        date={1980},
        ISSN={0097-3165,1096-0899},
     journal={J. Combin. Theory Ser. A},
      volume={28},
      number={2},
       pages={169\ndash 190},
         url={https://doi.org/10.1016/0097-3165(80)90083-7},
}

\bib{ContrerasMehtaSpan}{article}{
      author={Contreras, Ivan},
      author={Keller, Molly},
      author={Mehta, Rajan~Amit},
       title={Frobenius objects in the category of spans},
        date={2022},
        ISSN={0129-055X,1793-6659},
     journal={Rev. Math. Phys.},
      volume={34},
      number={10},
       pages={Paper No. 2250036, 34},
         url={https://doi.org/10.1142/S0129055X22500362},
}

\bib{CMS}{article}{
      author={Contreras, Ivan},
      author={Mehta, Rajan~Amit},
      author={Stern, Walker~H.},
       title={Frobenius and commutative pseudomonoids in the bicategory of spans},
        date={2025},
        ISSN={0393-0440,1879-1662},
     journal={J. Geom. Phys.},
      volume={207},
       pages={Paper No. 105309, 31},
         url={https://doi.org/10.1016/j.geomphys.2024.105309},
}

\bib{cooper-young}{incollection}{
      author={Cooper, Benjamin},
      author={Young, Matthew~B},
       title={Hall algebras via 2-{S}egal spaces},
        date={2026},
   booktitle={Higher {S}egal spaces and applications},
      series={Contemp. Math.},
      volume={838},
   publisher={Amer. Math. Soc., Providence, RI},
       pages={97\ndash 130},
         url={https://doi.org/10.1090/conm/838/16738},
}

\bib{DayStreet}{article}{
      author={Day, Brian},
      author={Street, Ross},
       title={Monoidal bicategories and {H}opf algebroids},
        date={1997},
        ISSN={0001-8708,1090-2082},
     journal={Adv. Math.},
      volume={129},
      number={1},
       pages={99\ndash 157},
         url={https://doi.org/10.1006/aima.1997.1649},
}

\bib{Dyckerhoff-Kapranov:Higher}{book}{
      author={Dyckerhoff, Tobias},
      author={Kapranov, Mikhail},
       title={Higher {S}egal spaces},
      series={Lecture Notes in Mathematics},
   publisher={Springer, Cham},
        date={2019},
      volume={2244},
         url={https://doi.org/10.1007/978-3-030-27124-4},
}

\bib{FGKW:unital}{article}{
      author={Feller, Matthew},
      author={Garner, Richard},
      author={Kock, Joachim},
      author={Proulx, May~U.},
      author={Weber, Mark},
       title={Every 2-{S}egal space is unital},
        date={2021},
        ISSN={0219-1997},
     journal={Commun. Contemp. Math.},
      volume={23},
      number={2},
       pages={Paper No. 2050055, 6},
         url={https://doi.org/10.1142/S0219199720500558},
}

\bib{GKT1}{article}{
      author={G\'{a}lvez-Carrillo, Imma},
      author={Kock, Joachim},
      author={Tonks, Andrew},
       title={Decomposition spaces, incidence algebras and {M}\"{o}bius inversion {I}: {B}asic theory},
        date={2018},
        ISSN={0001-8708},
     journal={Adv. Math.},
      volume={331},
       pages={952\ndash 1015},
         url={https://doi.org/10.1016/j.aim.2018.03.016},
}

\bib{GKT_comb}{incollection}{
      author={G\'{a}lvez-Carrillo, Imma},
      author={Kock, Joachim},
      author={Tonks, Andrew},
       title={Decomposition spaces in combinatorics},
        date={2026},
   booktitle={Higher {S}egal spaces and applications},
      series={Contemp. Math.},
      volume={838},
   publisher={Amer. Math. Soc., Providence, RI},
       pages={255\ndash 331},
         url={https://doi.org/10.1090/conm/838/16742},
}

\bib{goerss-jardine}{book}{
      author={Goerss, Paul~G.},
      author={Jardine, John~F.},
       title={Simplicial homotopy theory},
      series={Progress in Mathematics},
   publisher={Birkh\"auser Verlag, Basel},
        date={1999},
      volume={174},
        ISBN={3-7643-6064-X},
         url={https://doi.org/10.1007/978-3-0348-8707-6},
}

\bib{grothendieck:FGA3}{article}{
      author={Grothendieck, Alexander},
       title={Techniques de construction et th{\'e}oremes d'existence en g{\'e}om{\'e}trie alg{\'e}brique. {III}. {P}r{\'e}schemas quotients},
        date={1961},
     journal={S{\'e}minaire Bourbaki},
      volume={6},
      number={212},
       pages={99\ndash 118},
}

\bib{gurski:coherence}{book}{
      author={Gurski, Nick},
       title={Coherence in three-dimensional category theory},
      series={Cambridge Tracts in Mathematics},
   publisher={Cambridge University Press, Cambridge},
        date={2013},
      volume={201},
        ISBN={978-1-107-03489-1},
         url={https://doi.org/10.1017/CBO9781139542333},
}

\bib{hermida-mateus1}{article}{
      author={Hermida, Claudio},
      author={Mateus, Paulo},
       title={Paracategories. {I}. {I}nternal paracategories and saturated partial algebras},
        date={2003},
        ISSN={0304-3975,1879-2294},
     journal={Theoret. Comput. Sci.},
      volume={309},
      number={1-3},
       pages={125\ndash 156},
         url={https://doi.org/10.1016/S0304-3975(03)00135-X},
}

\bib{hermida-mateus2}{article}{
      author={Hermida, Claudio},
      author={Mateus, Paulo},
       title={Paracategories. {II}. {A}djunctions, fibrations and examples from probabilistic automata theory},
        date={2004},
        ISSN={0304-3975,1879-2294},
     journal={Theoret. Comput. Sci.},
      volume={311},
      number={1-3},
       pages={71\ndash 103},
         url={https://doi.org/10.1016/S0304-3975(03)00317-7},
}

\bib{joyal-street}{article}{
      author={Joyal, Andr\'e},
      author={Street, Ross},
       title={The geometry of tensor calculus. {I}},
        date={1991},
        ISSN={0001-8708,1090-2082},
     journal={Adv. Math.},
      volume={88},
      number={1},
       pages={55\ndash 112},
         url={https://doi.org/10.1016/0001-8708(91)90003-P},
}

\bib{lack:coherent}{article}{
      author={Lack, Stephen},
       title={A coherent approach to pseudomonads},
        date={2000},
        ISSN={0001-8708,1090-2082},
     journal={Adv. Math.},
      volume={152},
      number={2},
       pages={179\ndash 202},
         url={https://doi.org/10.1006/aima.1999.1881},
}

\bib{leinster}{article}{
      author={Leinster, Tom},
       title={Notions of {M}\"obius inversion},
        date={2012},
        ISSN={1370-1444,2034-1970},
     journal={Bull. Belg. Math. Soc. Simon Stevin},
      volume={19},
      number={5},
       pages={911\ndash 935},
}

\bib{MM:finset}{article}{
      author={Marx, Sophia~E},
      author={Mehta, Rajan~Amit},
       title={{Coherent span-valued 2D TQFTs}},
        date={2025},
     journal={arXiv preprint arXiv:2505.22832},
}

\bib{mateus-sernadas}{incollection}{
      author={Mateus, Paulo},
      author={Sernadas, Am\'ilcar},
      author={Sernadas, Cristina},
       title={Precategories for combining probabilistic automata},
        date={1999},
   booktitle={C{TCS} '99: {C}onference on {C}ategory {T}heory and {C}omputer {S}cience ({E}dinburgh)},
      series={Electron. Notes Theor. Comput. Sci.},
      volume={29},
   publisher={Elsevier Sci. B. V., Amsterdam},
       pages={Paper No. 29016, 18},
}

\bib{Morton}{article}{
      author={Morton, Jeffrey~Colin},
       title={Two-vector spaces and groupoids},
        date={2011},
     journal={Applied Categorical Structures},
      volume={19},
      number={4},
       pages={659\ndash 707},
}

\bib{rezk:homotopy}{article}{
      author={Rezk, Charles},
       title={A model for the homotopy theory of homotopy theory},
        date={2001},
        ISSN={0002-9947,1088-6850},
     journal={Trans. Amer. Math. Soc.},
      volume={353},
      number={3},
       pages={973\ndash 1007},
         url={https://doi.org/10.1090/S0002-9947-00-02653-2},
}

\bib{rota:comb1}{article}{
      author={Rota, Gian-Carlo},
       title={On the foundations of combinatorial theory. {I}. {T}heory of {M}\"obius functions},
        date={1964},
     journal={Z. Wahrscheinlichkeitstheorie und Verw. Gebiete},
      volume={2},
       pages={340\ndash 368},
         url={https://doi.org/10.1007/BF00531932},
}

\bib{segal:classifying}{article}{
      author={Segal, Graeme},
       title={Classifying spaces and spectral sequences},
        date={1968},
        ISSN={0073-8301,1618-1913},
     journal={Inst. Hautes \'Etudes Sci. Publ. Math.},
      number={34},
       pages={105\ndash 112},
         url={http://www.numdam.org/item?id=PMIHES_1968__34__105_0},
}

\bib{Segal:Conf}{article}{
      author={Segal, Graeme},
       title={Configuration-spaces and iterated loop-spaces},
        date={1973},
        ISSN={0020-9910},
     journal={Invent. Math.},
      volume={21},
       pages={213\ndash 221},
         url={https://doi.org/10.1007/BF01390197},
}

\bib{stay}{article}{
      author={Stay, Michael},
       title={Compact closed bicategories},
        date={2016},
     journal={Theory Appl. Categ.},
      volume={31},
       pages={Paper No. 26, 755\ndash 798},
}

\bib{Stern:span}{article}{
      author={Stern, Walker~H.},
       title={2-{S}egal objects and algebras in spans},
        date={2021},
        ISSN={2193-8407},
     journal={J. Homotopy Relat. Struct.},
      volume={16},
      number={2},
       pages={297\ndash 361},
         url={https://doi.org/10.1007/s40062-021-00282-8},
}

\bib{stern:perspectives}{incollection}{
      author={Stern, Walker~H.},
       title={Perspectives on the 2-segal conditions},
        date={2026},
   booktitle={Higher {S}egal spaces and applications},
      series={Contemp. Math.},
      volume={838},
   publisher={Amer. Math. Soc., Providence, RI},
       pages={1\ndash 30},
         url={https://doi.org/10.1090/conm/838/16735},
}

\bib{Verdon}{article}{
      author={Verdon, Dominic},
       title={Coherence for braided and symmetric pseudomonoids},
        date={2017},
     journal={arXiv preprint arXiv:1705.09354},
}

\end{biblist}
\end{bibdiv}

\end{document}